\newif\ifpreprint
\preprinttrue
\ifpreprint 
\documentclass[twocolumn,natbib]{article}
\else 
\documentclass[twocolumn,natbib]{svjour3}
\fi 
\usepackage{amsmath,amssymb}
\usepackage{color}
\usepackage[affil-it]{authblk}
\usepackage{algorithmic,algorithm}
\usepackage{graphicx}
\usepackage{hyphenat}
\usepackage{stmaryrd} 
\usepackage{enumitem}
\usepackage[T1]{fontenc}
\ifpreprint 
\usepackage[cm]{fullpage}
\fi 

\renewcommand{\vec}[1]{\boldsymbol{#1}}
\renewcommand{\matrix}[1]{\boldsymbol{#1}}

\newcommand{\Rearth}{R_{\operatorname{earth}}}
\newcommand{\Hatmos}{H_{\operatorname{atmos}}}

\newcommand{\order}[1]{\mathcal{O}(#1)}

\newcommand{\figdir}{./}

\newcommand{\CUSPARSE}{CUSPARSE}
\newcommand{\LAPACK}{LAPACK}
\newcommand{\CUBLAS}{CUBLAS}

\title{Matrix-free GPU implementation of a preconditioned conjugate gradient solver for anisotropic elliptic PDEs}
\ifpreprint 
\author[*,1]{Eike M\"uller}
\author[2]{Xu Guo}
\author[1]{Robert Scheichl} 
\author[2,3]{Sinan Shi}
\affil[1]{Department of Mathematical Sciences, University of Bath, Bath BA2 7AY, United Kingdom}
\affil[2]{Edinburgh Parallel Computing Centre (EPCC), The University of Edinburgh, James Clerk Maxwell Building, Mayfield Road, Edinburgh EH9 3JZ, United Kingdom}
\affil[3]{Current address: OT-Med, Europ\^{o}le M\'{e}diterran\'{e}en de l'Arbois, B\^{a}timent du Cerege, BP 80 13545 Aix-en-Provence Cedex 4, France}
\affil[*]{Email: \texttt{e.mueller@bath.ac.uk}}

\else 
\author{Eike M\"uller\thanks{Email: \texttt{e.mueller@bath.ac.uk}, Tel.: +44 1225 38 5633, Fax: +44 1225 38 6492} \and Xu Guo \and Robert Scheichl \and Sinan Shi}
\institute{Eike M\"uller%
\and%
Robert Scheichl%
\at%
Department of Mathematical Sciences, University of Bath, Bath BA2 7AY, United Kingdom%
\and%
Xu Guo%
\at%
Edinburgh Parallel Computing Centre (EPCC), The University of Edinburgh, James Clerk Maxwell Building, Mayfield Road, Edinburgh EH9 3JZ, United Kingdom%
\and%
Sinan Shi%
\at%
Edinburgh Parallel Computing Centre%
\and%
Current address: OT-Med, Europ\^{o}le M\'{e}diterran\'{e}en de l'Arbois, 
B\^{a}timent du Cerege, BP 80 13545 Aix-en-Provence Cedex 4, France }
\fi 
\begin{document}
\ifpreprint 
\twocolumn[
\begin{@twocolumnfalse}
\fi 
\maketitle 
\abstract{
Many problems in geophysical and atmospheric modelling require the fast solution of elliptic partial differential equations (PDEs) in ``flat'' three dimensional geometries. In particular, an anisotropic elliptic PDE for the pressure correction has to be solved at every time step in the dynamical core of many numerical weather prediction (NWP) models, and equations of a very similar structure arise in global ocean models, subsurface flow simulations and gas and oil reservoir modelling. The elliptic solve is often the bottleneck of the forecast, and to meet operational requirements an algorithmically optimal method has to be used and implemented efficiently.
Graphics Processing Units (GPUs) have been shown to be highly efficient (both in terms of absolute performance and power consumption) for a wide range of applications in scientific computing, and recently iterative solvers have been parallelised on these architectures.
In this article we describe the GPU implementation and optimisation of a Preconditioned Conjugate Gradient (PCG) algorithm for the solution of a three dimensional anisotropic elliptic PDE for the pressure correction in NWP.
Our implementation exploits the strong vertical anisotropy of the elliptic operator in the construction of a suitable preconditioner.
As the algorithm is memory bound, performance can be improved significantly by reducing the amount of global memory access. We achieve this by using a matrix-free implementation which does not require explicit storage of the matrix and instead recalculates the local stencil. Global memory access can also be reduced by rewriting the PCG algorithm using loop fusion and we show that this further reduces the runtime on the GPU. We demonstrate the performance of our matrix-free GPU code by comparing it both to a sequential CPU implementation and to a matrix-explicit GPU code which uses existing CUDA libraries. The absolute performance of the algorithm for different problem sizes is quantified in terms of floating point throughput and global memory bandwidth.
\ifpreprint 
\else 
\fi 
} 
\ifpreprint 
\end{@twocolumnfalse}
\vspace{3ex}
]
\fi 
\section{Introduction}
Anisotropic elliptic PDEs arise in many areas of geophysical and atmospheric modelling, which are often characterised by ``flat'' geometries: the horizontal extent of the domain of interest is much larger than its vertical size. This is the case for global weather- and climate prediction models. As the height of the atmosphere is significantly smaller than the radius of the earth, the horizontal resolution is of the order of 10-25 kilometers but the vertical grid spacing is several tens or hundreds of metres. Similar ranges of scales are encountered in models for simulating global ocean currents.
Due to the layered structure of geological formations, oil and gas reservoir simulations and subsurface flow models of aquifers are also typically carried out in ``flat'' domains. After discretisation the cells of the computational grid are very flat and the resulting matrix stencil is highly anisotropic, i.e. the coupling in the vertical direction exceeds the horizontal coupling by several orders of magnitude. 
To achieve optimal performance, it is important to exploit the strong anisotropy of the system when choosing an appropriate computational grid and an efficient solver. 

In this work we focus on the elliptic PDE for the pressure correction arising in the dynamical core of numerical weather- and climate- prediction models. In many forecast models semi-implicit semi\hyp Lagrangian time stepping introduced in \cite{Kwizak1971} and \cite{Robert1981} is used to advance the atmospheric fields forward in time. In contrast to explicit time stepping this method has a larger stability region and allows for longer model time steps without compromising the accuracy of the large scale dynamics, which can reduce the overall model runtime. However, if this approach is used to solve the fully compressible non-hydrostatic Euler equations, a three dimensional PDE for the pressure correction has to be solved at every time step as discussed for example in \cite{Smolarkiewiecz1997a,Thomas1997,Skamarock1997,Davies05,Melvin2010,Wood2013}, which often forms the computationally most expensive proportion of the model runtime.

Algorithmically the most efficient solvers for large elliptic PDEs are suitably preconditioned Krylov subspace\hyp\ or multigrid methods (see e.g. \cite{Hestenes1952,Briggs2000,Trottenberg2001,Saad03}). 
The strong anisotropy in the vertical direction can be exploited to construct an efficient preconditioner (or multigrid smoother) based on vertical line relaxation as discussed in \cite{Thomas1997,Skamarock1997}. In an related context \cite{Marshall1997} and \cite{Fringer2006} describe how the equations of ocean flows can be discretised on a tensor product grid which is unstructured in the horizontal but consists of regular columns in the vertical direction. Again the strong anisotropy is used in the construction of an efficient preconditioner of the iterative solver. In a similar fashion anisotropic elliptic PDEs arise in fully implicit methods for gas- and oil reservoir modelling. A ``supercoarsening'' multigrid algorithm for solving elliptic PDEs encountered in multiphase flow in porous media is described by \cite{Lacroix2003}: while the full three dimensional equation is solved on the finest grid, any vertical variations are averaged out on the coarser multigrid levels by collapsing vertical columns to a single layer.

The exact hardware on which forecast models will be implemented in the future is currently unknown, and it is important to explore novel chip architectures in addition to standard CPUs. Graphics Processing Units (GPUs) are fast and power- efficient computing devices and significant speedups relative to standard CPU implementations have been achieved in the past for iterative solvers for elliptic PDE, as described in \cite{Bolz03,Menon2007,Carvalho2010,Knittel2010,Dehnavi2011,Knibbe2011,deJong2012,Helfenstein2012,Reguly2012,Li2013}. 

While modern multicore CPUs contain several tens of cores and have a peak floating point performance of \mbox{$\order{10-100}$} GFLOPs, GPUs have several hundreds to thousands of cores and applications have to make efficient use of the massively parallel SIMD architecture and limited cache size per thread. The nVidia M2090 Fermi GPU, on which this work was carried out, has a peak performance of 1.331/0.665 TFLOP/s in single and double precision respectively and a global memory bandwidth of 177GByte/s. By dividing the peak FLOP rate by the memory bandwidth on the GPU, one can deduce that the number of computations per floating point variable read from memory is around 30, so computations are essentially ``free'' and the performance is limited by the speed with which data can be read from global memory and how efficiently it can be kept in cache. 

In this article we describe a matrix\hyp free GPU implementation of a preconditioned Conjugate Gradient (PCG) solver tailored towards the solution of an\-iso\-tro\-pic PDEs with a tensor-product structure. 
The most compute intensive components of the iterative solver are the evaluation of a large sparse matrix-vector product (SpMV) and the inversion of a block-tridiagonal matrix. Both kernels were ported to the GPU and the memory access pattern and thread layout were adapted to increase data throughput. 
In our implementation the matrix is not stored explicitly but recalculated in every grid cell, which reduces access to global memory compared to the matrix-explicit code. However the stencil we use is more complicated than the simple Poisson stencil on a regular grid used in previous studies (see \cite{Menon2007,Knittel2010}). Due to the tensor product structure of the equation and the computational grid, the matrix can be written as the product of a one dimensional vertical discretisation, which is the same for every column and only needs to be calculated and stored once, and a horizontal stencil. As the number of vertical levels is very large in atmospheric applications, and the horizontal coupling only needs to be calculated once for each vertical column, this creates only a small overhead. The vertical discretisation requires the storage of four vectors of length $n_z$, where $n_z$ is the number of vertical levels. In total we store $4\times n_z$ values to parametrise the matrix.
With the exception of \cite{Menon2007,Knittel2010}, all implementations discussed in the literature that we are aware of, store the matrix explicitly, which requires the storage of $7\times n_{\operatorname{horiz}}\times n_z$ matrix entries, where $n_{\operatorname{horiz}}$ is the number of horizontal grid cells. This can have a negative impact on the performance on bandwidth limited architectures as it requires reading the matrix stencil from global memory in addition to the field vectors.
While the specific implementation described in this article is based on a three dimensional grid which can be written as the tensor product of regular horizontal grid and a graded vertical grid, the method we present can also be applied to unstructured horizontal grids.

As the sparse matrix-vector product and preconditioner solve are highly efficient in our GPU implementation, other parts of the main CG such as level 1 BLAS vector updates and scalar products start to account for a significant proportion of the runtime, even if they are implemented using optimised GPU libraries such as \CUBLAS. We find that to achieve further performance increases it is not sufficient to optimise the kernels in isolation, but rather several components of the main CG iteration need to be considered together, as has been suggested in \cite{Dehnavi2011}. In particular we find that the number of memory references can be reduced further by fusing the loops over the computational grid in the main kernels with the BLAS operations and this can lead to an additional performance gain of around $30\%$. For the entire PCG algorithm we are able to obtain a total speedup of a factor $60\times$ for single precision arithmetic on a nVidia Fermi M2090 card relative to one core of an Intel Xeon Sandybridge E5-2620 CPU. For double precision arithmetic the speedup was slightly smaller with $48\times$. This includes time for setting up the discretisation and copying data between host and device. To study the performance of our matrix\hyp free GPU code we compared it to an implementation which stores the matrix explicitly in the compressed sparse row storage (CSR) format using the \CUSPARSE\ and \CUBLAS\ libraries. Our matrix\hyp free code is significantly faster than the implementation based on CRS data structures which does not exploit the regular structure of the problem. We quantified the absolute performance in terms of floating point operations per second (FLOPs) and global memory bandwidth for different problem sizes, where the latter is the more relevant measure for the performance of a memory bound algorithm. The optimised matrix-free code achieves a bandwidth of around $25\%-50\%$ of the theoretical peak value, which is a sizeable proportion but shows that theoretically there is still potential for additional improvements which could lead to a further speedup of a factor $2\times-4\times$. An idea of how this could be achieved by increasing the granularity of the algorithm is discussed below. The floating point performance is 70-80 GFLOPs for single precision and 40-50 GFLOPs for double precision, corresponding to around $5-8\%$ of the theoretical peak value.
\\[2ex]
\paragraph{Overview.} This article is organised as follows: previous GPU implementations of iterative solvers for PDEs are reviewed in section \ref{sec:Review}. In section \ref{sec:ModelEquation} the model equation and its discretisation is described in detail with particular emphasis on the tensor-product structure of the grid and the elliptic operator. Preconditioned Krylov-subspace solvers and the matrix-free and interleaved form of the PCG algorithm for solving the model equation are presented in section \ref{sec:KrylovSolvers}. A general overview over the GPU architecture and the CUDA programming- and execution model can be found in section \ref{sec:GPUs}, and our CUDA implementation of the PCG solver is described in section \ref{sec:Implementation}.
Performance measurements are discussed in section \ref{sec:NumericalExperiments} where we also present comparisons to a matrix-explicit implementation and quantify the absolute performance. Our conclusions and a discussion of planned further work can be found in section \ref{sec:Conclusion}. For reference appendix \ref{sec:Algorithms} contains the explicit form of the two most important kernels of our optimised algorithm.
\section{Previous work}\label{sec:Review}
The GPU implementation of Krylov-subspace solvers and in particular of the Preconditioned Conjugate Gradient algorithm has been studied extensively in the literature, both for more general sparse matrices and for matrices arising from the discretisation of elliptic PDEs. As far as we are aware, all implementations discussed in the literature (with the exception of \cite{Knittel2010} and \cite{Menon2007}) are based on matrix-explicit representations.
While some of the authors study Poisson- or sign-positive Helmholtz equations, none of the problems studied in the literature show the strong anisotropy which characterises the elliptic operator we consider in this work, and hence the preconditioners investigated in the literature will not be optimal in our case.

While the speedups presented in the following review depend on the problem and on the hardware used at a particular time and should only be used as an indicator for achievable performance gains, almost all GPU implementations are significantly faster than the corresponding CPU versions with speedups of $20\times-40\times$ relative to the sequential code.

Some early work is presented in \cite{Bolz03} where both a conjugate gradient and a multigrid solver are implemented for solving the sign positive Helmholtz equation $-\Delta u+\sigma u=g$ arising from an implicit time discretisation of the incompressible Navier-Stokes equations on a regular two dimensional grid. However, for both solvers the matrix is stored explicitly.

The compressed sparse row storage format (CSR) is a very popular and general format which has been used in a variety of recent GPU implementations of Krylov subspace algorithms.
A PCG solver for the same sign-positive Helmholtz equation as in \cite{Bolz03} is described in \cite{Helfenstein2012} for two and three dimensions. For the approximate inverse SSOR preconditioner, which requires an additional matrix multiplication, a socket-to-socket speedup of more than $8\times$ is reported for the best implementation of the PCG algorithm on an nVidia Tesla T10 card relative to the unpreconditioned CG algorithm on an Intel Xeon Quad-Core 2.66 GHz CPU. Although in contrast to \cite{Bolz03} a three dimensional system is solved, the elliptic operator considered is fully isotropic in both cases.
In \cite{Dehnavi2011} a modified version of the Conjugate Gradient algorithm is used for solving a set of general matrices from the University of Florida sparse matrix collection described in \cite{Davis2011}. A simple Jacobi preconditioner is used and the performance of the solver is optimised by using the prefetch CSR sparse matrix-vector multiplication in \cite{Dehnavi2010} and fusing kernels in the main PCG loop. As described in section \ref{sec:PCGilvd} below we use a similar technique for fusing different kernels in our implementation to improve the performance of the code. Together with some other improvements the authors of \cite{Dehnavi2011} report that this led to a significant speedup compared to a GPU implementation using the ``Row per warp'' sparse matrix-vector multiplication described in \cite{Bell2009}. One of the problems studied in \cite{Dehnavi2011} is the ``thermal2'' matrix which arises from an FEM discretisation of the stationary heat equation $\partial_x(k\partial_x T)+\partial_y(k \partial_y T)=0$ on an unstructured two dimensional grid. For this problem a speedup of $41\times$ relative to a single core of a Intel Core2 2.4GHz could be achieved on both nVidia GT8800 and GTX280 GPUs.
The GPU implementation of a Krylov subspace solver for the (sign-indefinite) two dimensional Helmholtz equation is described in \cite{Knibbe2011}. A shifted Laplace multigrid preconditioner is used to reduce the number of iterations and a speedup of around $30\times$ could be achieved on an nVidia GeForce 9800 GTX/9800 GTX+ GPU relative to the sequential implementation on one core of an AMD Phenom 9850 CPU.

Other sparse matrix storage formats have also been used to implement iterative solvers on GPUs. An implementation of a CG solver for the two dimensional PDE arising from the implicit time discretisation of the heat equation is described in \cite{Michels2011}. The ELLPACK-R data format, which is more suitable for structured problems, was used for storing the matrix and a speedup of a factor $26\times$ could be achieved for a two dimensional problem of size $2048\times 2048$ on an nVidia GeForce GTX 480 card, relative to the sequential implementation on an Intel Core i7 860 CPU with 2.80GHz. Although the structure of the five point nearest-neighbour stencil arising from a finite-difference approximation of the Poisson equation is similar to the stencil we use in our discretisation, in contrast to our problem the elliptic PDE solved in \cite{Michels2011} is two dimensional and fully isotropic. The GPU implementation of preconditioned GMRES and Conjugate Gradient solvers for a range of problems and preconditioners has been studied in \cite{Li2013} and the performance for different sparse matrix storage formats is compared. While in some of the implementations only the sparse-matrix vector product is carried out on the device and the preconditioner is executed on the host, preconditioners that are easier to parallelise are also ported to the GPU. However, the authors find that for a simple block-Jacobi preconditioner implemented on the GPU the number of iterations is very large. This should be compared to our implementation: for the strongly anisotropic elliptic PDE we consider the blocks have a direct physical interpretation as they described the strong vertical coupling within one column, which is much larger than the coupling between different blocks. As a result, the simple block-diagonal preconditioner proved to be very efficient in our numerical tests. In \cite{deJong2012} the GPU implementation of a Preconditioned Conjugate Gradient solver for both a two dimensional Poisson problem and the elliptic equation arising in the Variational Boussinesq Model (VBM) is described. A Repeated Red Black (RRB) preconditioner is used, but an incomplete Poisson preconditioner with diagonal scaling is also considered. As for the sparse approximate inverse used in \cite{Helfenstein2012}, the incomplete Poisson preconditioner can be reduced to an additional sparse matrix product. The sparse matrix-vector product is implemented by storing the local five-point stencil at each gridpoint. For the RRB preconditioner a speedup of around $40\times$ could be achieved for the Poisson test problem on an nVidia GeForce GTX 580 card, relative to the sequential implementation on an Intel Xeon W3520 CPU. However, the costs for memory allocation and setup of the preconditioner matrix take up around a third of the total runtime. In contrast in our matrix free implementation only a small amount of data has to be copied between host and device and the matrix setup costs are negligible. For realistic problems the speedup reported in \cite{deJong2012} is $20\times-30\times$ for the RRB preconditioner and $5\times-20\times$ for the incomplete Poisson preconditioner.

As far as we are aware, the only matrix-free implementation of a CG solver discussed in the literature are \cite{Menon2007,Knittel2010}. In \cite{Menon2007} both the homogenous Poisson equation and the Navier Stokes equation are solved on an unstructured mesh by implementing matrix-free gradient and divergence operations. On an nVidia a speedup of around $3\times$ was achieved on an nVidia 6600GT card relative to an AMD Athlon 64 CPU. Note, however, that the implementation is based on the low-level graphics API and the hardware used in the study is quite dated by current standards. The GPU implementation of a matrix-free PCG solver for the homogenous Poisson equation in three dimensions is also described in \cite{Knittel2010}.

Due to its significance in many scientific applications and in particular iterative solvers, the performance of sparse matrix-vector multiplications on its own has been studied extensively in the literature:
Various sparse matrix formats are described in \cite{Bell2009} and their performance for a sparse matrix-vector multiplication is compared for both structured and unstructured matrices. While CSR is the most general format and can be applied to matrices with widely varying row sizes, the best performance for structured matrices arising from the discretisation of PDEs is obtained with the DIA and ELLPACK formats. However, in \cite{Reguly2012} an efficient parameter dependent implementation of sparse matrix-vector multiplication based on the CSR format on cache based GPUs is described. Cache usage and performance can be improved significantly by varying the number of threads processing each row, the thread block size and number of rows processed by one cooperating thread group. The authors find that by tuning the parameters heuristically, the performance of a Conjugate Gradient solver for a structured problem arising from the finite element discretisation of a simple elliptic Poisson problem using CSR storage is comparable to the corresponding implementation using the ELLPACK format. Both matrix-explicit versions are beaten by an implementation which stores a small local matrix on each element and assembles the global stiffness matrix on-the-fly in each matrix-vector product as described in \cite{Markall2010,Cantwell2011}. For other work on improved matrix-explicit implementations of the sparse matrix-vector product see the review and references cited in \cite{Reguly2012}. 

The number of iterations can often be reduced significantly by using multigrid methods, and recently both geometric and algebraic multigrid solvers have been ported to GPUs, see e.g. \cite{Goodnight2005,Geveler2011,Brannick2013}. The extension of our PCG solver to a geometric multigrid solver for an\-iso\-tro\-pic problems based on the tensor product idea in \cite{BoermHiptmair1999} is discussed in \cite{Mueller2013} and we are currently working on a GPU implementation of the same matrix-free geometric multigrid solver.

We finally remark that multiple-GPU implementations of iterative solvers have been described in the literature (see e.g. \cite{Cevahir2009,Georgescu2010,Knittel2010,Verschoor2012}) and we are currently working on extending our algorithm to clusters of GPU.
\section{Model equation}\label{sec:ModelEquation}
Following \cite{Melvin2010,Wood2013} a model equation which reproduces the most important features of the full PDE for the pressure correction in a NWP model, has been derived in \cite{Mueller2013}. The derivation of the full pressure correction equation in atmospheric models can be found for example in 
\cite{Smolarkiewiecz1997a,Thomas1997,Skamarock1997,Davies05} and is described for ocean models in \cite{Marshall1997} and \cite{Fringer2006}. 
The model equation we use is a symmetric positive definite PDE and can be written in spherical coordinates as
\begin{equation}
  -\omega^2\left(\Delta_{2d}+\lambda^2\frac{1}{r^2}\frac{\partial}{\partial r}\left(r^2\frac{\partial}{\partial r}\right)\right)u + u = f
\label{eqn:ModelEquation}
\end{equation}
where $\Delta_{2d}$ denotes the two dimensional Laplacian on the unit sphere. 
For simplicity all length scales are measured in units of the earth radius $\Rearth$, and the equation is solved in a thin spherical shell $r\in \left[1,1+\Hatmos\right]$. We write $\Hatmos=D/\Rearth\ll 1$ where $D$ is the thickness of the atmosphere.
The model parameters $\omega^2$ and $\lambda^2$ depend on the model time step size, in particular, $\omega$ is proportional to the time step size.
In our numerical experiments the parameters were adjusted to their values for typical model resolutions in NWP with a Courant number of around $10$.
An important feature of the discretised PDE is the strong anisotropy in the vertical direction: the depth of the atmosphere is about two orders of magnitude smaller than the radius of the earth and consequently the horizontal grid spacing $\Delta x$ is significantly larger than the vertical grid spacing $\Delta z$. The strong grid aligned anisotropy in the vertical direction is given (approximately) by
\begin{equation}
  \gamma^2 = \left(\frac{\lambda\; \Delta x}{\Delta z}\right)^2,
  \label{eqn:Anisotropy}
\end{equation}
and as $H_{\operatorname{atmos}}\ll 1$ we have $\Delta z \ll \Delta x$, so $\gamma^2 \gg 1$. 
This property can be used to construct a simple but very efficient and parallelisable preconditioner for iterative solvers of this equation, which solves the vertical equation exactly but ignores the horizontal couplings.

The condition number of the preconditioned operator approaches a fixed value as the horizontal resolution increases and does not diverge as for the Poisson equation. To see this, note that to keep the Courant number constant, the time step size $\Delta t\propto\omega$ has to chosen to be proportional to the horizontal grid spacing, i.e. $\Delta t\propto \Delta x$. The largest eigenvalue of the preconditioned matrix is of the order $\omega^2/\Delta x^2\propto \Delta t^2/\Delta x^2$ and independent of $\Delta x$, whereas the smallest eigenvalue is 1 due to the presence of the zero order term in (\ref{eqn:ModelEquation}). 

After discretisation the elliptic operator in (\ref{eqn:ModelEquation}) can be written abstractly in tensor product form as the sum of three terms:
\begin{equation}
  L = D^{(h)}\otimes M^{(r)} + M^{(h)}\otimes D^{(r)} + \tilde{M}^{(h)}\otimes \tilde{M}^{(r)}
  \label{eqn:TensorProductOperator}
\end{equation}
Here $M^{(h)}$ and $\tilde{M}^{(h)}$ ($M^{(r)}$ and $\tilde{M}^{(r)}$) are horizontal (vertical) mass matrices which only contain couplings in the horizontal (vertical) direction. $D^{(h)}$ ($D^{(r)}$) are second order derivate operators which contain couplings in the horizontal (vertical) direction. 
\subsection{Discretisation}
To discretise (\ref{eqn:ModelEquation}) we use a finite volume scheme on a tensor-product grid, which consists of a (possibly unstructured but conforming) two dimensional grid on the unit sphere and a non-uniform (typically graded) one dimensional grid in the vertical direction. The model fields are defined as integrals in a grid cell given by the horizontal grid element $T$ and vertical index $k=0,\dots,n_z-1$ \begin{equation}
  \overline{u}_{k}^{(T)} \equiv \int_{r_k}^{r_{k+1}}\int_T u(r,\vec{\theta})\;r^2 dr\;d\vec{\theta}
\label{eqn:ColumnRepresentation}
\end{equation}
with $\vec{\theta}$ denoting horizontal coordinates on the unit sphere (throughout this work we use zero-based indexing as all our implementations are in the C programming language).
Independent of the horizontal discretisation we need to store a vector $\overline{\vec{u}}^{(T)}$ of length $n_z$ at each element $T$.
The vertical grid is defined by the grid points $r_k$ where $k=0,\dots,n_z$. The number of vertical grid cells is usually large, $n_z=\order{100}$.
In meteorological applications a graded vertical grid with smaller grid spacings near the ground is desirable and we use $r_k=1+(k/n_z)^2\cdot \Hatmos$.

Equation (\ref{eqn:ModelEquation}) is discretised using a finite volume scheme, and schematically it can be written for each horizontal grid element $T$ as
\begin{equation}
  (\matrix{A}\overline{\vec{u}})^{(T)} = \matrix{A}_T \overline{\vec{u}}^{(T)} + \sum_{T'\in\mathcal{N}(T)}\matrix{A}_{T,T'}\overline{\vec{u}}^{(T')} = \overline{\vec{f}}^{(T)}
\label{eqn:ModelEqnCubedSphere}
\end{equation}
where $\mathcal{N}(T)$ is the set of neighbours of $T$.
$\matrix{A}_T$ is a tridiagonal $n_z\times n_z$ matrix and $\matrix{A}_{T,T'}$ are diagonal $n_z\times n_z$ matrices for each neighbouring element $T'$. 
The matrices $\matrix{A}_{T,T'}$ correspond to the off-diagonal couplings in the horizontal derivative operator $D^{(h)}\otimes M^{(r)}$ in (\ref{eqn:TensorProductOperator}) and are given by the product
\begin{equation}
  \matrix{A}_{T,T'} = \alpha_{T,T'}\operatorname{diag}(\vec{d}) 
 =
  \alpha_{T,T'}
\begin{pmatrix}
    d_0 & & &\\
        & d_1 & &   \\
      &  & \ddots & \\
      & & & d_{n_z-1}
  \end{pmatrix}.
  \label{eqn:DiagonalMatrix}
\end{equation}
The coefficient $\alpha_{T,T'}$ is the ratio between the length $S_{T,T'}$ of the edge between the cells $T$ and $T'$ and the distance between the midpoints $\vec{r}_{T}$ and $\vec{r}_{T'}$ of these cells. We also define the diagonal coefficient $\alpha_{T}$ as the sum of the off-diagonal coefficients
\begin{equation}
  \alpha_{T} \equiv \sum_{T'\in\mathcal{N}(T)} \alpha_{T,T'}.
\end{equation}
The (symmetric) tridiagonal matrix $\matrix{A}_{T}$ can be split into a sum of three terms:
\begin{equation}
  \begin{aligned}
  \matrix{A}_{T}&=  \begin{pmatrix}
    \tilde{d}_0 & b_0 &  \\
    c_1 & \ddots & \ddots  \\
      & \ddots & \ddots & b_{n_z-2} \\
      & & c_{n_z-1} & \tilde{d}_{n_z-1}
  \end{pmatrix}\\
&=|T|\operatorname{diag}(\vec{a}) - \alpha_{T} \operatorname{diag}(\vec{d}) \\&\qquad+ |T|\operatorname{tridiag}\left(-(\vec{b}+\vec{c}),\vec{b},\vec{c}\right),
 \end{aligned}
  \label{eqn:TridiagonalMatrix}
\end{equation}
with $\tilde{d}_k = |T|(a_k-b_k-c_k)-\alpha_{ij}d_k$ where $|T|$ is the area of the horizontal grid element $T$. The first term corresponds to the product $\tilde{M}^{(h)}\otimes\tilde{M}^{(r)}$ in (\ref{eqn:TensorProductOperator}). The second term is the diagonal contribution of $D^{(h)}\otimes M^{(r)}$, and the third term corresponds to the vertical derivative term $M^{(h)}\otimes D^{(r)}$. 
While the finite volume discretisation described so far leads to a seven-point nearest-neighbour stencil on a regular grid, the same structure also arises for other stencil types which include couplings between grid cells which are not directly adjacent.

The vectors $\vec{a}$, $\vec{b}$, $\vec{c}$ and $\vec{d}$ do not depend on the grid cell $T$ and can be precomputed. On the other hand the quantities $|T|$ and $\alpha_{T,T'}$ only need to be calculated once per vertical column. These two important observations will be exploited in the efficient matrix-free implementation of the PCG solver described below. 
\subsection{Global matrix representation}
\begin{figure}
 \begin{center}
   \includegraphics[width=0.9\linewidth]{\figdir/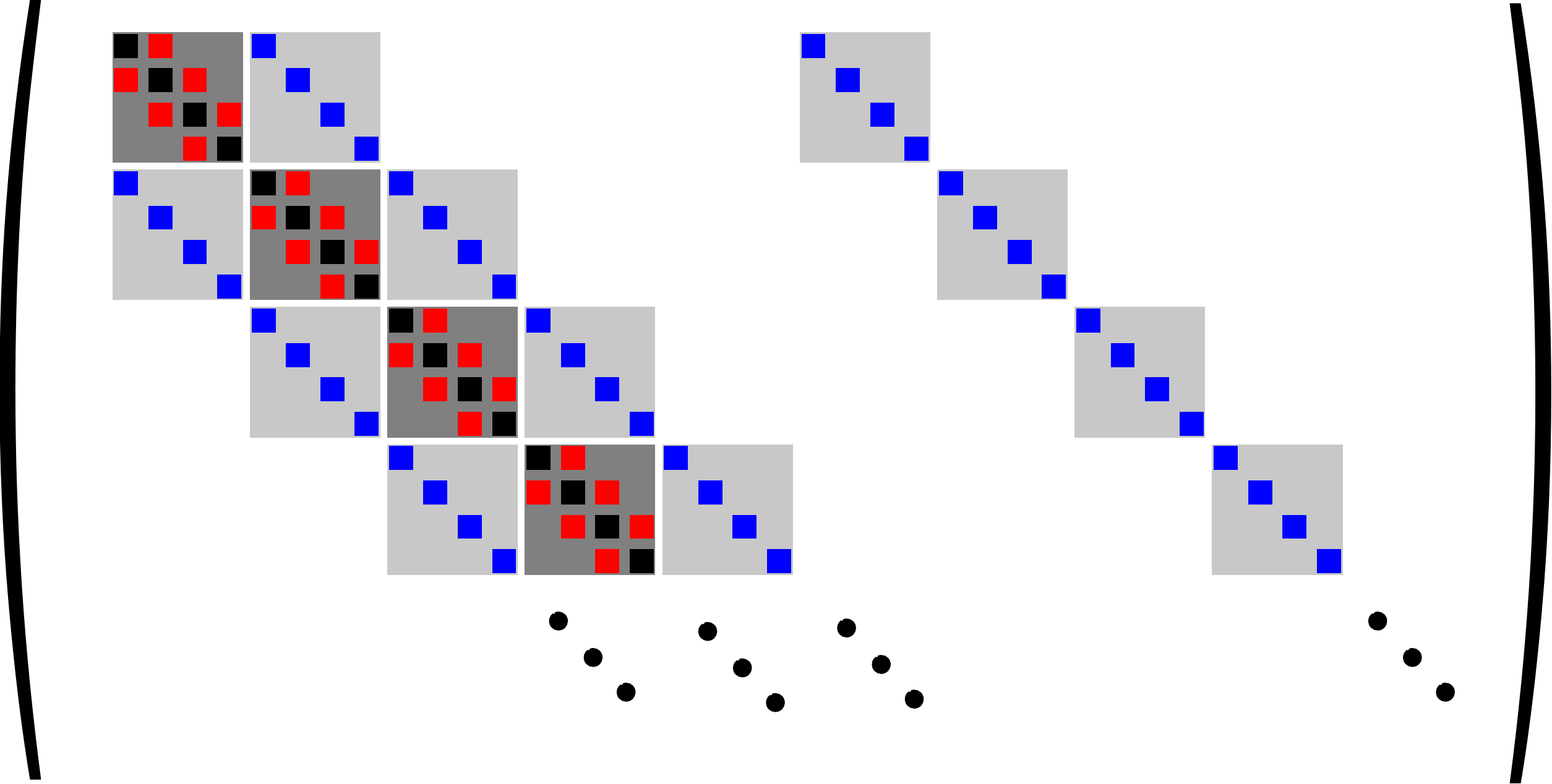}
   \caption{Structure of the matrix arising from the finite volume discretisation of (\ref{eqn:ModelEquation}) for $n_z=4$ vertical columns. Vertical couplings are shown in red and horizontal couplings in blue.}
   \label{fig:MatrixStructure}
 \end{center}
\end{figure}

Assuming an ordering of the horizontal degrees of freedom, which maps each cell $T$ to a linear index $\nu(T)\in 0,\dots,n_{\operatorname{horiz}}-1$, one can write the full 3d solution vector of length $n=n_{\operatorname{horiz}}\times n_z$ as
\begin{equation}
  \vec{u} = \left\{\overline{\vec{u}}^{(0)},\overline{\vec{u}}^{(1)},\dots,\overline{\vec{u}}^{(n_{\operatorname{horiz}})}\right\}
\end{equation}
with
\begin{equation}
  u_{\ell} = \overline{u}^{(T)}_{k} \qquad \text{where}\quad \ell = n_z\cdot \nu(T)+k.
\end{equation}
In this representation the discretised equation (\ref{eqn:ModelEqnCubedSphere}) can be written in the familiar form as 
\begin{equation}
  \matrix{A}\vec{u} = \vec{f}
  \label{eqn:MatrixVectorEquation}
\end{equation}
and the structure of the matrix $\matrix{A}$ is shown in Fig. \ref{fig:MatrixStructure}. Note that the matrix has a block structure, where each of the blocks corresponds to one combination $(T,T')$ of neighbouring elements of the horizontal grid.
Each of the gray blocks is of size $n_z\times n_z$, the dark gray blocks ($T'=T$) describe the diagonal terms and vertical coupling, whereas the light gray blocks ($T\neq T'\in\mathcal{N}(T)$) describe the horizontal coupling.
\begin{figure}
  \begin{center}
  \includegraphics[width=0.6\linewidth]{\figdir/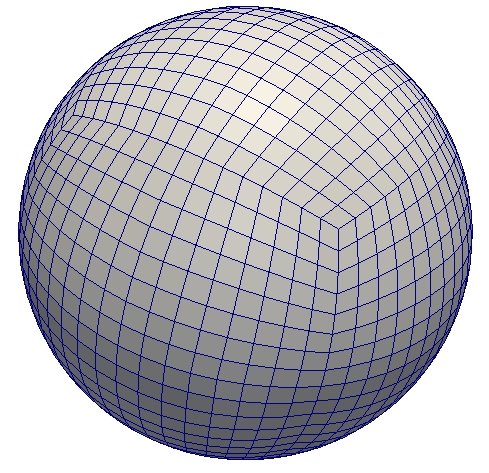}
  \caption{Cubed sphere grid. The model equation was discretised on one logically rectangular panel of this grid.}
  \label{fig:GridCubedSphere}
  \end{center}
\end{figure}
In the following we work on one panel of a cubed sphere grid with gnomonic projection (Fig. \ref{fig:GridCubedSphere}). This defines a mapping from $\tilde{\Omega}=[-1,1]\times[-1,1]$ to spherical coordinates on one-sixth of the entire sphere and each cell of the regular grid of size $n_{\operatorname{horiz}}=m\times m$ on $\tilde{\Omega}$ can be labelled with an index pair $(i,j) \in [0,m-1]\times [0,m-1]$. In this case we have
\begin{equation}
  u_{\ell} = \overline{u}^{(i,j)}_{k} \qquad \text{where}\quad \ell = n_z\left(m\cdot i + j\right)+k
\label{eqn:GlobalColumnReprConversion}
\end{equation}
and label each grid cell $T_{ij}$ by its horizontal indices. 
We stress, however, that there is no algorithmic problem in extending our approach to unstructured horizontal grids or to the entire spherical shell for example by ordering the horizontal grid cells along a space-filling curve.
\section{Iterative solvers for elliptic PDEs}\label{sec:KrylovSolvers}
Typically the number of degrees of freedom per atmospheric variable in current global forecast models is at the order of several 10 millions. For next generation forecast models global horizontal model resolutions of around one kilometre are envisaged, which will require the solution of PDEs with more than $10^{10}$ unknowns.
Clearly naive direct methods can not be used for the solution of equations of this size and spectral methods often require a regular underlying grid structure and are hard to parallelise. 

Krylov subspace methods (see e.g. \cite{Saad03}) are very efficient iterative algorithms for solving sparse linear systems, in particular if they are suitably preconditioned.
These methods construct the approximation $\vec{u}^{(k)}$ of the exact solution $\vec{u}$ of (\ref{eqn:MatrixVectorEquation}) in a $k$-dimensional Krylov-subspace
\begin{equation}
\mathcal{K}_k=\operatorname{span}\left\{\vec{r},\matrix{A}\vec{r},\matrix{A}^2\vec{r},\dots,\matrix{A}^{k-1}\vec{r}\right\} \subset \mathbb{R}^n,
\end{equation}
where $\vec{r}$ is the initial residual $\vec{r}=\vec{b}-\matrix{A}\vec{u}^{(0)}$. They are easy to parallelise as in addition to local operations such as sparse matrix-vector multiplications and \texttt{axpy}-like vector updates they require only a small number of global reductions.
\subsection{Preconditioned Conjugate Gradient algorithm}
The simplest Krylov subspace method, which can be applied for symmetric positive definite matrices, is the Conjugate Gradient (CG) algorithm introduced in \cite{Hestenes1952}. At each step the approximate solution vector $\vec{u}^{(k)}$ is updated by adding a correction proportional to the search direction $\vec{p}^{(k)}$. The search directions are chosen such that they are $\matrix{A}$- orthogonal for different $k$, $k'$: $\langle \vec{p}^{(k)},\matrix{A}\vec{p}^{(k')}\rangle = 0$ for $k\ne k'$. The convergence rate depends on the spectral properties of the matrix $\matrix{A}$, in particular on the condition number $\kappa$, which is the ratio between the largest and smallest eigenvalue, as derived in \cite{Saad03}.
For the system arising from the discretisation of the Poisson equation, $\kappa$ grows rapidly with the inverse grid spacing $1/h$ and it can be shown that asymptotically for $h\rightarrow 0$ the number of iterations required to reduce the error by a factor $\epsilon$ is
\begin{equation}
  k \propto \frac{\log \epsilon}{h}.\label{eqn:IterationsCG}
\end{equation}
For anisotropic systems $h$ is the smallest grid spacing in the problem, for the highly anisotropic problem in this work this is the vertical grid spacing $\Delta z \ll \Delta x$.
The dominant cost in each iteration is the sparse matrix-vector multiplication $\vec{y}\mapsfrom \matrix{A}\vec{x}$, which is of $\order{n}$ computational complexity. Hence the total cost of the algorithm is
\begin{equation}
  \operatorname{Cost}({\operatorname{CG}}) \propto \frac{n}{h}\log \epsilon.
  \label{eqn:CostKrylov}
\end{equation}
To solve non-symmetric systems, more general Krylov space methods such as GMRES, BCG, BiCGStab and GCR can be used. Although the number of sparse matrix-vector products and intermediate vectors which need to be stored changes between different Krylov subspace algorithms, their general structure is very similar and in this work we focus on the Conjugate Gradient algorithm for simplicity.

An equivalent version of the linear system in (\ref{eqn:MatrixVectorEquation}) can be obtained by multiplication with the inverse of the matrix $\matrix{M}$,
\begin{equation}
  \matrix{M}^{-1}\matrix{A}\vec{u} = \matrix{M}^{-1}\vec{f}.
\end{equation}
This is generally referred to as left preconditioning. $\matrix{M}$ is a matrix which should be easy to invert (i.e. it should be easy to solve the system $\matrix{M}\vec{x}=\vec{y}$) and as similar to $\matrix{A}$ as possible, such that the preconditioned matrix $\matrix{M}^{-1}\matrix{A}$ is well conditioned.
Usually these two requirements are mutually exclusive and a tradeoff between them has to be found. Initial profiling of the CPU code revealed that the sparse matrix-vector multiplication and preconditioner solve account for the largest proportion of the runtime ($80-90\%$).
In addition to these operations each iteration requires three \texttt{axpy}-like vector updates, one \texttt{scal}-operation and two scalar products (\texttt{dot}). An additional norm calculation (\texttt{nrm}) is required for the evaluation of the stopping criterion.
The number of floating points operations and memory references for the individual level 1 BLAS operations is given in Tab. \ref{tab:BLASoperations}, the total number of FLOPS per grid cell for all BLAS operations is $13$ and the number of memory references is $16$, and hence this part of the algorithm is clearly memory bound.
\begin{table}
  \begin{center}
  \begin{tabular}{|l|rr|}
    \hline
    operation & FLOPs & MEM \\\hline\hline
    \texttt{scal} & $1$ & $2$\\
    \texttt{axpy} & $2$ & $3$\\
    \texttt{dot}  & $2$ & $2$\\
    \texttt{nrm2} & $2$ & $1$\\
    \textbf{total (PCG)} & $\boldsymbol{13}$ & $\boldsymbol{16}$\\
    \hline
  \end{tabular}
  \caption{Number of floating point operations and memory references per grid cell for different level 1 BLAS operations. The last row shows the total number of FLOPs and memory references for all BLAS operations in the PCG algorithm.
}
  \label{tab:BLASoperations}
  \end{center}
\end{table}

The strong anisotropy of the discretised PDE can be used to construct an efficient preconditioner: if the small horizontal couplings are ignored, the matrix $\matrix{A}$ shown in Fig. \ref{fig:MatrixStructure} is block-diagonal with each block corresponding to one vertical column.
Each of the tridiagonal blocks can be inverted independently using the Thomas algorithm written down explicitly in \cite{Press2007}. 
More efficient block-SOR preconditioners can be used as well, and require the inversion of the same block-diagonal matrix plus an additional sparse matrix-vector product.

This approach has been applied very successfully for the pressure solver in the dynamical core of several numerical weather- and climate prediction models, see \cite{Smolarkiewicz1994,Skamarock1997,Davies05,Piotrowski2011}. In particular the authors of \cite{Thomas1997} show the efficiency of a simple 1d line relaxation in comparison to other preconditioners such as 2d ADI or a three dimensional pointwise SOR iteration.
The good weak and strong scaling on up to 65536 cores of the HECToR Cray supercomputer and more than $10^{10}$ degrees of freedom is demonstrated for the model equation (\ref{eqn:ModelEquation}) in \cite{Mueller2013}.

As discussed in section \ref{sec:ModelEquation}, the condition number of the preconditioned elliptic operator approaches a fixed value for physical choices of the parameter $\omega$ in (\ref{eqn:ModelEquation}) as the horizontal grid spacing tends to zero. Hence for these parameters Krylov subspace solvers for this equation are algorithmically stable in the sense that the number of iterations does not diverge as the horizontal resolution increases. However, as shown in \cite{Mueller2013} the number of iterations and total solution time can be reduced significantly by using (geometric) multigrid solvers.
\subsection{Matrix-free implementation}\label{sec:MatrixFreeImplementation}
Neither the matrix $\matrix{A}$ nor the preconditioner matrix $\matrix{M}$ need to be stored explicitly in the algorithm, it is sufficient to evaluate the sparse matrix-vector product $\vec{y}\mapsfrom \matrix{A}\vec{x}$ and to solve the equation $\matrix{M}\vec{x}=\vec{y}$ for $\vec{x}$. For matrices arising from the discretisation of PDEs, such as the one discussed in section \ref{sec:ModelEquation}, the local matrix stencil only couples each grid cell to its neighbours. As memory access is significantly more expensive than floating point operations on GPUs, it will be beneficial to recalculate the stencil whenever it is needed in the sparse matrix-vector product or preconditioner solve. 
In each horizontal grid cell $(i,j)$ the diagonal matrices $\matrix{A}_{T_{ij},T_{i'j'}}$ in (\ref{eqn:DiagonalMatrix}) (with $T_{i'j'}\in\mathcal{N}(T_{ij})$) and the tridiagonal matrix $\matrix{A}_{T_{ij}}$ in (\ref{eqn:TridiagonalMatrix}) are calculated from the precomputed vectors $\vec{a}$, $\vec{b}$, $\vec{c}$ and $\vec{d}$. For efficiency, we instead store the vectors $\vec{a}'$, $\vec{b}'$, $\vec{c}'$ and $\vec{d}$ with 
\begin{xalignat}{3}
  a'_k &= a_k/d_k,&
  b'_k &= b_k/d_k,&
  c'_k &= c_k/d_k,
\end{xalignat}
as then the number of floating point operations in the sparse matrix-vector multiplication can be further reduced. To reconstruct the matrix stencil the coefficients $|T_{ij}|$ and $\alpha_{i'j'}$ are also needed (we write $\alpha_{ij}\equiv\alpha_{T_{ij}}$ and $\alpha_{i'j'}\equiv\alpha_{T_{ij},T_{i'j'}}$ for $T_{i'j'}\in \mathcal{N}(T_{ij})$ for simplicity). While these can be given by relatively complicated algebraic expressions on complex geometries, such as the spherical grid considered in this article, they only need to be calculated once in each vertical column and for large $n_z$ this will only lead to a small overhead.

On the regular horizontal grid which we use in our implementation the sparse matrix-vector product $\vec{y}\mapsfrom \matrix{A}\vec{x}$ can then be calculated in each grid cell $(i,j,k)$ as
\begin{equation}
\begin{aligned}
  y_{ijk} &\mapsfrom \Big[\left((a'_k-b'_k-c'_k)\cdot |T_{ij}| - \alpha_{ij}\right)\cdot x_{ijk}\\
  & \quad+\; |T_{ij}|\cdot b'_k\cdot x_{i,j,k+1}+ |T_{ij}|\cdot c'_k \cdot x_{i,j,k-1}\\
  & \quad+\; \alpha_{i+1,j}\cdot x_{i+1,j,k}
  +\alpha_{i-1,j}\cdot x_{i-1,j,k}\\
  &\quad+\alpha_{i,j+1}\cdot x_{i,j+1,k}
  +\alpha_{i,j-1}\cdot x_{i,j-1,k}
  \Big]\cdot d_k
\end{aligned}
\end{equation}
(with possible modifications at the boundary of the domain).
For each $(i,j,k)$ this requires 20 floating point operations and 12 memory references (7 loads for $\vec{x}$, 1 store for $\vec{y}$ and 4 loads for $a'_k$, $b'_k$, $c'_k$ and $d_k$). However, as the vectors $\vec{a}'$, $\vec{b}'$, $\vec{c}'$ and $\vec{d}$ do not change from column to column these are likely to remain in cache, reducing the number of memory references to 8. Depending on the innermost loop two of the values of $\vec{x}$ (namely the ones belonging to the same vertical column which are needed at the next vertical level, i.e. $x_{ijk}$ and $x_{ij,k+1}$) will most likely stay in cache, which reduces the number of memory references further to 6.
In contrast 14 floating point operations and 22 memory references are necessary if the matrix $\matrix{A}$ is stored in compressed sparse row storage (CSR) format. Again, the actual number of memory references is likely to be smaller as some of the variables are cached. However, in this case caching is not possible for the matrix entries which vary from one three dimensional grid cell to the next.

For the construction of the preconditioner we drop the second term in 
(\ref{eqn:ModelEqnCubedSphere}), which couples different vertical columns, and write $\matrix{A}_{T_{ij}}\overline{\vec{u}}^{(i,j)} = \overline{\vec{b}}^{(i,j)}$ (using (\ref{eqn:GlobalColumnReprConversion}) to implicitly convert between the global vector representation and the column based representation in (\ref{eqn:ColumnRepresentation})). For each column $(i,j)$ the tridiagonal system $\matrix{A}_{T_{ij}}\overline{\vec{x}}^{(i,j)}=\overline{\vec{y}}^{(i,j)}$ defined by
\begin{equation}
\begin{aligned}
  |T_{ij}|\cdot d_k\cdot \Big[&\quad\overline{x}^{(i,j)}_{k-1}\cdot c'_k\\
 &+\overline{x}^{(i,j)}_{k}\cdot \left((a'_k-b'_k-c'_k)-\tilde{\alpha}_{ij}\right)\\
 &+\overline{x}^{(i,j)}_{k+1}\cdot b'_k\Big] = \overline{y}^{(i,j)}_{k}
\end{aligned}
\label{eqn:MatrixFreeTridiagonalSystem}
\end{equation}
with $\tilde{\alpha}_{ij}\equiv \alpha_{ij}/|T_{ij}|$ needs to be solved for $\overline{\vec{x}}^{(i,j)}$. This can be done efficiently in $\order{n_z}$ time using the Thomas algorithm, which is essentially Gaussian elimination applied to a tridiagonal system. 
The forward iteration requires 8 memory references at each step (loading the auxilliary vector $\phi_{k-1}$ and $y^{(i,j)}_{k-1}$ and storing $\phi_k$ and $y^{(i,j)}_k$ plus four loads for $a'_k$, $b'_k$, $c'_k$ and $d_k$). In each step of the backward iteration 4 memory references are needed. Hence the total number of memory references is $12$. Again, some of the data might be kept in cache. If only the vectors $\vec{a}'$, $\vec{b}'$, $\vec{c}'$ and $\vec{d}$ are cached the number of memory references reduces to 8. If in addition any data in the same vertical column can be kept in cache, the number of memory references reduces further to 5. As there are no dependencies in the horizontal direction, we can parallelise in this direction, i.e. the tridiagonal solve in each vertical column can be carried out independently.

An additional advantage of the matrix free method is the fact that there  are no matrix setup costs; the cost for precomputing the vectors $\vec{a}'$, $\vec{b}'$, $\vec{c}'$, $\vec{d}$ and copying them to the device is neglible as these vectors only have length $n_z$.
\subsection{Interleaved PCG algorithm}\label{sec:PCGilvd}
Field vectors are accessed in different components of the PCG algorithm, for example the residual $\vec{r}$ is needed in the preconditioner solve, in the update of the residual and the calculation of the residual norm.
In the standard implementation of the algorithm these operations are carried out in separate loops over the grid, which increases the number of memory references as data can not be kept in cache. However, the main iteration of the PCG algorithm can be rewritten such that it only consists of two loops over the grid, each of which contains either the sparse matrix-vector multiplication or the tridiagonal solve and a number of BLAS operations. Similar loop fusion for the a GPU implementation of the PCG algorithm presented in \cite{Chronopoulos1989} has been discussed in \cite{Dehnavi2011}.

The main iteration of this modified algorithm is shown in Algorithm \ref{alg:InterleavedPCG}.
\begin{algorithm}
  \caption{Interleaved PCG loop}
  \label{alg:InterleavedPCG}
  \begin{algorithmic}[1]
    \FOR{$k=1,\operatorname{maxiter}$}
    \STATE{Interleaved preconditioner kernel: Calculate\\
$\vec{r}\mapsfrom \vec{r}-\alpha \vec{q}$,
$\vec{z}=\matrix{M}^{-1}\vec{r}$,
$||\vec{r}||\mapsfrom \sqrt{\langle \vec{r},\vec{r}\rangle}$,
$\kappa\mapsfrom  \langle \vec{r},\vec{z}\rangle$\\
      in a single iteration over the grid.}
    \STATE{\textbf{if} ($||\vec{r}||/||\vec{r}_0|| < \epsilon$ or $||\vec{r}||<\tau)$ \textbf{then} Exit}
    \STATE{$\beta\mapsfrom \kappa/\kappa_{old}$, $\kappa_{old}\mapsfrom \kappa$}
    \STATE{Interleaved SpMV kernel: Calculate\\
        $\vec{u}\mapsfrom \vec{u}+\alpha \vec{p}$,
        $\vec{p}\mapsfrom \vec{z}+\beta \vec{p}$,
        $\vec{q} \mapsfrom \matrix{A}\vec{z}+\beta \vec{q}$,
        $\sigma\mapsfrom\langle \vec{p},\vec{q}\rangle$\\
      in a single iteration over the grid.}
    \STATE{$\alpha \mapsfrom \kappa_{old}/\sigma$}
   \ENDFOR
\end{algorithmic}
\end{algorithm}
The kernels are written down explicitly for the matrix-free implementation in Appendix \ref{sec:Algorithms}. The number of floating point operations and memory references for the matrix-free PCG algorithm and for its interleaved version are shown in Tab. \ref{tab:PCGOperations}. For the memory references we give three different values corresponding to the following assumptions: (1) no data is kept in cache, (2) only the vectors $\vec{a'}$, $\vec{b'}$, $\vec{c'}$ and $\vec{d}$ are cached and (3) in addition data in the same vertical column is cached.

While the algorithm is still memory bound on the GPU, the number of memory references is reduced in the interleaved implementation, in particular if the cache can be used efficiently.
\begin{table}
 \begin{center}
  \begin{tabular}{|l||ll|lll|}
  \hline 
  algo-& operation & FLOPs & \multicolumn{3}{|c|}{Memory references}\\
  rithm          &           &       & no & matrix & columns\\
            &           &       & cache & cached & cached\\\hline\hline
  PCG & SpMV  & 20 & 12  & 8 & 6 \\
      & prec  & 13 & 12  & 8 & 5 \\
      & BLAS  & 13 & 16  & 16 & 16 \\
      & \textbf{total} & \textbf{46} & \textbf{40} & \textbf{32} & \textbf{27}\\
  \hline
  Inter- & SpMV & 28 & 17 & 13 & 11 \\
  leaved & prec & 19 & 16 & 12 & 9 \\
  PCG & \textbf{total} & \textbf{47} & \textbf{33} &\textbf{25} &\textbf{20}\\
  \hline
  \end{tabular}
  \caption{Number of floating point operations and memory references per iteration and per grid point for different components of the PCG and the interleaved PCG algorithm. Memory references are shown without caching, with caching the vectors $\vec{a'}$, $\vec{b'}$, $\vec{c'}$ and $\vec{d}$ only and assuming that all degrees of freedom in a vertical column are cached as well.
}\label{tab:PCGOperations}
 \end{center}
\end{table}
\section{Graphics Processing Units}\label{sec:GPUs}
In contrast to CPUs, on which most transistors are used for advanced execution control units and cache hierarchies, GPUs have a very large number (several hundreds to thousands) of lightweight compute cores which support SIMD parallelism for simple compute kernels and are ideally suited for floating point intensive calculations on regular data. The clock speed of each individual core is smaller than on the CPU, which improves the power efficiency.

The cores in the GPU are grouped into a number of streaming multiprocessors (SMs). Data can be stored in global on-chip GPU memory and in addition, each SM has a smaller and faster shared memory. Each compute core can also access an even smaller local memory in addition to a set of registers. Constants can be stored in fast constant memory. Modern GPUs, such as the Fermi architecture (see \cite{FermiWhitepaper2009}), also have a hardware managed L1/L2 cache hierachy.
Data transfer between host and device memory has to be managed explicitly by the user. As typically the PCIe bus has a small bandwidth (at the order of ten GB/s), this is expensive and data should be calculated and kept on the GPU as long as possible.
\subsection{CUDA programming model}
The CUDA programming model described in the CUDA programming guide (\cite{CUDA2012}) provides an extension of the C language. When writing code for a GPU, the most compute intensive subroutines are parallelised by isolating them in simple kernels. On the GPU these kernels are executed by light\-weight threads which are run on the compute cores of the SMs. Threads are grouped into blocks in a one-, two- or three dimensional grid, which execute independently on one SM. Synchronisation and data exchange is only possible between threads in the same block. In each block the threads are arranged into a grid, i.e. each thread is uniquely identified by its thread index and block index.
Threads are scheduled in groups of size 32, called warps. The threads in each warp execute one common instruction at a time; if the execution path diverges, for example due to an \texttt{if}-statement in the kernel, both branches are executed. This is known as thread divergence and should be avoided whenever possible. 

As our solver algorithm is memory bound, performance can be increased by minimising latency and increasing global memory bandwidth.
If a warp stalls as it is waiting for data from memory, it is paused and the next warp which is ready for execution is launched by the warp scheduler. In contrast to threads on a multicore CPU, the GPU scheduler is designed for launching and switching threads with minimal overheads. Thus, given the number of threads is large enough, memory latency can be hidden. To achieve this the GPU usually has to be oversubscribed, i.e. the number of threads launched at a single time should exceed the number of compute cores.

For optimal efficiency, memory access for all threads in one half-warp should be coalesced. Global memory access is processed in segments of 128 bytes (which is the size of one cache line) and all threads in a half-warp should access the smallest possible number of different segments. For example, if each of the 16 threads reads a double precision (8 byte) floating point number from global memory, this will require the transfer of one segment if these numbers are stored consecutively, but it will require 16 separate memory transfers and thus incur a big penalty if the numbers are more than 128 byte apart in global memory.
\section{CUDA implementation of the PCG algorithm}\label{sec:Implementation}
Both the standard PCG method and its interleaved version were implemented in C and CUDA-C. In the standard version the sparse matrix-vector multiplication and preconditioner were implemented as kernels and the \CUBLAS\ library was used for the level 1 BLAS operations.

To minimise memory transfers between host and device, data is kept on the device inside the entire PCG loop and all operations are carried out on the GPU. Host-device data transfers are only necessary for copying the initial solution and the right hand side to GPU memory before the CG iteration and for copying the final solution back to the host at the end. In addition, the four vectors $\vec{a}'$, $\vec{b}'$, $\vec{c}'$ and $\vec{d}$ describing the vertical discretisation need to be copied to the device once before the main PCG iteration. However, as each vector only has length $n_z$, the amount of data tranferred is negligible in comparison to the size of the initial solution and right hand side vector. In particular it is significantly less than for matrix-explicit implementations. 
\subsection{Domain decomposition and data layout}\label{sec:MemoryLayout}
Due to the inherently sequential nature of the Thomas algorithm, parallelisation in the vertical direction is not possible for the preconditioner.
Instead the code is parallelised by assigning each vertical column to one thread and organising these into threadblocks of size $B=B_x\times B_y$, each of which is launched on one streaming multiprocessor. Parallelisation only in the horizontal direction is common in atmospheric models where data dependencies in the vertical direction are introduced by physical processes such as radiative transfer. To achieve a good occupancy the number of threads per block should be large, in particular latency hiding can only occur for $B\gg32$.

In the code all three dimensional field vectors, such as the solution vector $\vec{u}$, are represented as one dimensional arrays of size $n=m\times m\times n_z$ which are stored contiguously in memory.  To access the entry $u_{ijk}$ the three dimensional index $(i,j,k)\in[0,m-1]\times [0,m-1] \times [0,n_z-1]$ (where $k$ is the vertical index) has to be mapped to a linear index $\ell\in\{0,\dots,n-1\}$. To achieve optimal cache usage on the CPU, vertical columns are stored consecutively in memory in the C code, which can be achieved by the mapping
\begin{equation}
  \ell^{\operatorname{(C)}} = n_z\cdot (m\cdot i+j)+k\label{eqn:MemoryMappingCPU}
\end{equation}
already introduced in (\ref{eqn:GlobalColumnReprConversion}). However, on a GPU, the data layout in (\ref{eqn:MemoryMappingCPU}) would lead to a significant amount of uncoalesced memory access as the number of vertical levels is large. For a typical $n_z$ of 128, consecutive threads will access data which is $128\times \operatorname{sizeof}(float)$ bytes apart in memory. Although this problem can be mitigated by use of the L1 cache (or manual prefetching into shared memory at the beginning of each kernel), on a GPU the L1 cache is shared between a large number of threads, which severly limits the cache size per thread.
As each thread processes an entire vertical column, efficient caching would only possible for relatively small horizontal block sizes $B$: in the Thomas algorithm two vectors of length $n_z$ need to be stored per column in addition to the four vectors $\vec{a}'$, $\vec{b}'$, $\vec{c}'$ and $\vec{d}$ which describe the vertical discretisation. The total amount of L1 cache on the Fermi architecture is limited to 48kB, hence for single precision floating point numbers one would have
\begin{equation}
  48\operatorname{kB} \ge (2B + 4)n_z\times \operatorname{sizeof}(float)
\end{equation}
which would limit the number of threads per block to 44 for single precision and 22 for double precision calculations if we assume that $n_z=128$.

\begin{figure}
  \begin{center}
    \includegraphics[width=0.8\linewidth]{\figdir/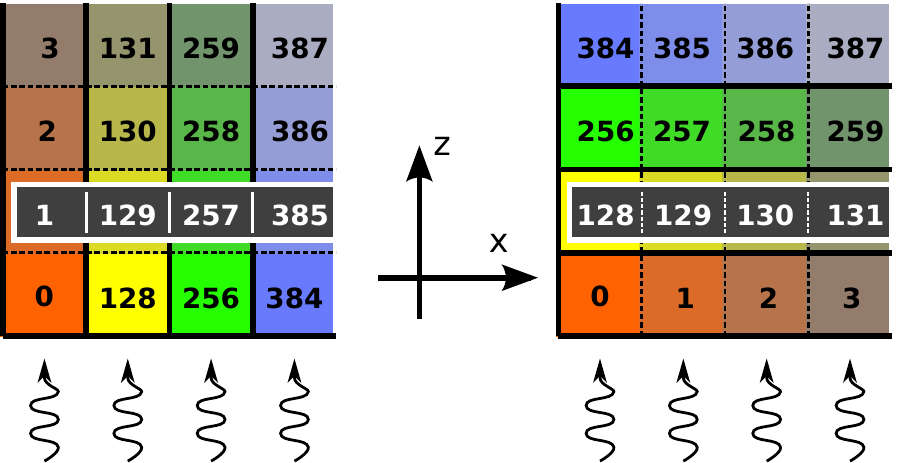}
    \caption{Memory layout and data access pattern. Data read by all threads in a warp is shown with a gray background. If the vertically continuous ordering of the degrees of freedom in (\ref{eqn:MemoryMappingCPU}) is used, memory access is not coalesced between the threads in one warp (left). The ordering degrees of freedom  in (\ref{eqn:MemoryMappingCUDAC}) leads to coalesced memory access between all threads in a warp (right).}
    \label{fig:CoalescedMemoryAccess}
  \end{center}
\end{figure}

An alternative approach is to change the storage format of the fields. Instead of (\ref{eqn:MemoryMappingCPU}) we use 
\begin{equation}
  \ell^{\operatorname{(CUDA-C)}} = m\cdot (n_z\cdot j+k)+i\label{eqn:MemoryMappingCUDAC}
\end{equation}
in the matrix-free CUDA-C code, i.e. the first horizontal index runs fastest. This ensures that, provided the horizontal block size $B_x$ is larger than 16, memory access for all threads in a half-warp is coalesced, as illustrated in Fig. \ref{fig:CoalescedMemoryAccess}.
In our numerical tests we found that the blocksize $B_x=64$, $B_y=2$ gave good results, in particular the global load efficiency was almost $100\%$ for the interleaved preconditioner and larger than $88\%$ for the SpMV kernel. For this blocksize the data processed by one block is too large to fit into cache. However, we find that the L1 cache hit rate ranges between $33\%$ for the interleaved preconditioner and $56\%$ for the interleaved SpMV kernel. The total global memory bandwidth is around $25\%-50\%$ of the peak value for both kernels. To reduce the number of cache misses further, a more fine grained parallelisation would have to be used to reduce the data volume per thread and ensure that data used by each thread block can be kept in shared memory.

The approach we are currently exploring (but which is not used in the implementation described in this article) is to use a different solver for the tridiagonal system in the vertical direction. 
The substructuring method discussed in \cite{Toselli2005} splits the $n_z\times n_z$ problem into $B_z$ smaller tridiagonal systems of size $\approx (n_z/B_z) \times (n_z/B_z)$, which can be solved independently by different threads. To obtain the global solution on the entire vector of length $n_z$, a global tridiagonal system of size $B_z\times B_z$ for the interface points has to be solved before the solutions of the small subsystems can be combined to the total solution. 
\subsection{Matrix-explicit implementation}
In addition to the matrix-free implementation described in section \ref{sec:MatrixFreeImplementation} we also wrote a version of the code based on an explicit representation of the matrix. Because of its popularity in the literature on sparse matrix-vector products and Krylov-subspace solvers (see the discussion in section \ref{sec:Review}) we chose the CSR format. The $n\times n$ tridiagonal matrix used in the preconditioner was stored as a set of three vectors of length $n$. The \verb!cusparse{S|D}csrmv()! function from the \CUSPARSE\ library was used for the sparse matrix-vector multiplication and the subroutine \verb!cusparse{S|D}gtsv()! for the tridiagonal solve in the preconditioner.

Setup of the matrix on the host and copying the CSR representation to the device would create additional costs as host-device data transfers are expensive. For a fair comparison the matrices were set up on the device instead. Our numerical tests show that in this case the setup costs only form a small part of the total runtime. However, this might not be the case in other applications where matrix has to be constructed on the CPU.

In addition we wrote a CPU version using our own hand-written CSR matrix-vector product and the \LAPACK\ routines \texttt{GTTRF}, \texttt{GTTRS} for the tridiagonal solver.
\section{Numerical experiments}\label{sec:NumericalExperiments}
\subsection{Hardware and compilers}
All runs were carried out on the GPU node of the \texttt{aquila} cluster in Bath. The node contains an Intel Xeon E5-2620 Sandybridge CPU with a clockspeed of 2.00GHz and an nVidia Fermi M2090 GPU. The theoretical peak performance of one core of the Sandybridge CPU without AVX extensions is 8.0GFLOPs (4 floating point operations per cycle $\times$ 2.0 GHz). The M2090 GPU contains 512 cores running at a clockspeed of 1.3GHz which are organised into 16 streaming multiprocessors with 32 cores each, as described in the Fermi Architecture Whitepaper (\cite{FermiWhitepaper2009}).
The total size of global GPU memory is 6GB, and each SM has 64kB of on-chip memory which can be split between shared memory and the L1 cache as 48kB/16kB or 16kB/48kB. In our implementation we used API calls to choose the optimal partitioning. In addition the Fermi architecture has 768kB of global L2 cache.
The theoretical peak performance is quoted as 1.331TFLOPs for single precision and 0.665 TFLOPs for double precision while the peak bandwidth for access to global memory is quoted as 177GB/sec. Dividing the peak performance by the peak bandwidth implies that around 30 floating point operations can be carried out for each single/double precision variable read from global memory. In practise the actual number of floating point operations per accessed variable will of course be different due to latency effects and as required data might already be in L1 cache. However, this again stresses the importance of optimising memory throughput to achieve good performance.

In contrast, on the CPU the theoretical peak bandwidth is 41.6GB/s, and the number of floating point operations per value loaded from memory is 0.76 for single- and 1.54 for double precision.

The nVidia \texttt{nvcc} compiler (release 5.0, V0.2.1221) was used for compiling the CUDA code and we used version 4.4.6 of the gnu C compiler for complilation of the CPU code. To achieve the best possible performance of the matrix-explicit CPU code, optimised BLAS and LAPACK libraries based on the OpenBLAS implementation (see \cite{OpenBLAS}) were used. AVX extensions were disabled on the CPU.

\subsection{Results}\label{sec:Results}
We first study the performance and speedup of the individual components of the PCG solver and of the entire algorithm for a fixed problem of size $256\times 256\times128$ with $\omega^2=6.71\cdot10^{-4}$ and $\lambda^2=3.32\cdot10^{-2}$, which is a typical set of parameters in meteorological applications. On the GPU we used a two-dimensional block layout and each block has a size of $B_x\times B_y = 64\times 2$ threads. We found that this gives good results and varying the block size did not increase the performance.
\subsection{Matrix-vector multiplication and preconditioner}
In Tab. \ref{tab:ResultsApplyPrecMatrixExplicit} the times for a single sparse matrix-vector multiplication and preconditioner solve are shown for the matrix-explicit method. These times do not include any costs for setting up the matrix or for transferring data between host and device, as this is only required once for each PCG solve, which consists of a large number of sparse matrix-vector multiplications and preconditioner applications. The speedups shown in this table are relative to the sequential CPU implementation. 
\begin{table}
 \begin{center}
  \begin{tabular}{|l|cc|c|}
  \hline
    \textbf{matrix-explicit} & \multicolumn{2}{|c|}{time per call} & speedup \\
    kernel         &  C & CUDA & C vs. CUDA \\
\hline\hline
    \textit{single precision} &&&\\
    SpMV         & 170.6 & 6.91 & $25\times$ \\
    preconditioner & 205.3 & 12.50 & $16\times$ \\
  \hline
    \textit{double precision} &&&\\
    SpMV         & 182.2 & 10.91 & $17\times$ \\
    preconditioner & 249.7 & 23.20 & $11\times$ \\
  \hline
  \end{tabular}\\[1ex]
  \caption{Measured times and speedups for one sparse matrix-vector multiplication (SpMV) $\vec{y}\mapsfrom \matrix{A}\vec{x}$ and one preconditioner solve $\vec{x}\mapsfrom \matrix{M}^{-1}\vec{y}$ on the CPU and GPU using the matrix-explicit implementation. All times are given in milliseconds. The speedup of the CUDA-C code relative to the sequential CPU implementation is shown in the last column for each case.}
  \label{tab:ResultsApplyPrecMatrixExplicit}
 \end{center}
\end{table}
Assuming that the CPU code can be parallelised perfectly, the socket-to-socket speedup is a factor $6$ smaller as the CPU contains six processor cores. 

The corresponding times for the matrix-free implementation are shown in Tab. \ref{tab:ResultsApplyPrecMatrixFree}, where we also show the speedup of the matrix-free CUDA-C code relative to the matrix-explicit CUDA-C code.
\begin{table}
 \begin{center}
  \begin{tabular}{|l|cc|cc|}
  \hline
    \textbf{matrix-free} & \multicolumn{2}{|c|}{time per call} & \multicolumn{2}{|c|}{speedup} \\
                   &   &        & C vs. & mat.-free \\
    kernel         & C & CUDA & CUDA & vs. CSR \\
\hline\hline
    \textit{single precision} &&&&\\
    SpMV         & 78.5 & 0.75 & $105\times$ & $9.2\times$\\
    preconditioner & 252.6 & 2.40 & $105\times$ & $5.2\times$\\
    interlvd. SpMV & 129.9 & 2.16 & $60\times$ & ---\\
    interlvd. prec. & 253.1 & 3.34 & $76\times$ & ---\\
\hline
    \textit{double precision} &&&&\\
    SpMV         & 80.7 & 1.41 & $57\times$ & $7.7\times$\\
    preconditioner & 350.4 & 3.75 & $93\times$ & $6.2\times$\\
    interlvd. SpMV & 132.3 & 3.86 & $34\times$ & ---\\
    interlvd. prec. & 351.3 & 4.86 & $72\times$ & ---\\
  \hline
  \end{tabular}
  \caption{Measured times and speedups for one sparse matrix-vector multiplication (SpMV) $\vec{y}\mapsfrom \matrix{A}\vec{x}$ and preconditioner solve $\vec{y}\mapsfrom \matrix{M}^{-1}\vec{x}$ on the CPU and GPU using the matrix-free implementation. All times are given in milliseconds. The last two columns show the speedup of the matrix-free CUDA-C code relative to the corresponding sequential CPU implementation and relative to the matrix-explicit CUDA-C version.}
  \label{tab:ResultsApplyPrecMatrixFree}
 \end{center}
\end{table}
On the CPU the matrix-free sparse matrix-vector multiplication is more than twice as fast as the matrix-explicit implementation. However, the matrix-explicit preconditioner is $25\%-30\%$ faster than the corresponding matrix-free version. On the GPU the matrix-free code is significantly faster than the matrix-explicit version, both for the sparse matrix-vector multiplication where the speedup is $9.2\times$ for single precision and $7.7\times$ for double precision. For the preconditioner the speedup is slightly smaller with $5.2\times$ and $6.2\times$ for single- and double precision respectively.
The speedup of both standalone matrix-free GPU kernels relative to the sequential CPU code is more than $100\times$ in single precision. In double precision the speedup is $93\times$ for the preconditioner and $57\times$ for the sparse matrix-vector multiplication. While the corresponding speedups for the interleaved kernels are smaller, their performance has to be judged in the context of the full PCG loop, which for the standalone kernels also contains several level 1 BLAS operations.

We observe that for the matrix-free standalone SpMV kernel the double precision implementation takes about twice as long as the single precision version. This suggests that the implementation is bandwidth limited and less affected by cache efficiency, as the double precision version requires transferring twice as much data from global memory than the single precision implementation. This interpretation is also corroborated by the relatively high cache hit rate for this kernel reported in Tab. \ref{tab:ResultsPerformancePrecMatrixFree}. The cache hit rate is smaller for the interleaved sparse matrix-vector multiplication and the preconditioner kernels, all of which show a smaller increase in the runtime between single and double precision.
A similar drop of performance by nearly a factor of two when going from single to double precision can be observed for the matrix-explicit kernels, see Tab. \ref{tab:ResultsApplyPrecMatrixExplicit}.
\subsection{PCG algorithm}\label{sec:ResultsPCG}
We now analyse the performance of the entire PCG solver and break down the time spent in a single iteration of the algorithm.
\subsubsection{Total solution time}
We measured the time required to carry out in 100 PCG iterations, which is sufficient to reduce the residual by five orders or magnitude. Our measurements include the time for the matrix setup and data transfer between host and device. The results are listed for three different implementations of the algorithm in Tab. \ref{tab:TotalTimes} where we also calculated the speedup of the matrix-free CUDA-C code both relative to the C code and relative to the matrix-explicit GPU code. The matrix-free interleaved algorithm gives the best performance, with a speedup of a factor of $60\times$ (single precision) and $48\times$ (double precision) relative to the C code on the CPU. It outperforms the matrix-explicit GPU code by more than a factor of four.
On the CPU the interleaved algorithm is only slightly faster than standard PCG for single precision and even slightly slower for double precision. This is because in the sequential implementation most of the time ($80\%-90\%$) is spent in the sparse matrix-vector multiplication and preconditioner kernels, so fusing the kernels can only give a speedup of no more than $10\%-20\%$. This is different on the GPU, as will be discussed in section \ref{sec:TimePerIteration}.

As expected the cost for setting up the vertical discretisation matrix (calculation of $\vec{a}'$, $\vec{b}'$, $\vec{c}'$ and $\vec{d}$) and copying it to the device turned out to be negligible ($\ll 1ms$).
For the matrix-explicit code the matrix setup time only accounts for a small proportion of the runtime; on the GPU the matrix setup time is $6\%$ and $4\%$ of the total solution time in single- and double precision respectively. Although for the matrix-free interleaved code host-device data transfer of the solution and right hand side vector takes up only a small part of the runtime (about $8\%$ both in single- and double precision), this is not true any longer if a smaller number of iterations is used for example to only reduce the residual by three orders of magnitude. We also found that on the CPU the total solution time can be reduced by a factor of around four if the Krylov- subspace solver is replaced by a geometric solver multigrid, as is discussed in \cite{Mueller2013}. Again, this would increase the relative importance of the host-device memory transfer.
\begin{table*}
 \begin{center}
  \begin{tabular}{|l|cc|c|cc|cc|}
  \hline
                 & \multicolumn{2}{|c|}{matrix and} & 
                    \multicolumn{1}{|c|}{data} &
                    \multicolumn{2}{|c|}{total time} & 
                    \multicolumn{2}{|c|}{speedup}\\ 
                    & \multicolumn{2}{|c|}{preconditioner setup}& 
                      transfer &
                      &  &
                     C vs. & matrix-free\\
    implementation & C & CUDA & 
                     CUDA &
                     C & CUDA &
                     CUDA & vs. CSR\\
\hline\hline
    \textit{single precision} &&&&&&&\\
    matrix-explicit & 0.55 & 0.15 & 0.047 & 43.72 & 2.54 & $17\times $ & ---\\
    matrix-free & --- & --- & 0.047 & 37.00 & 0.87 & $43\times $ & $2.9\times$\\
    matrix-free interleaved & --- & --- & 0.047 & 37.39 &  0.62 & $60\times $ & $4.1\times$\\
  \hline
    \textit{double precision} &&&&&&&\\
    matrix-explicit & 0.72 & 0.16 & 0.073 & 50.50 & 4.41 & $11\times $ & ---\\
    matrix-free & --- & --- & 0.073 & 49.37 & 1.48 & $33\times $ & $3.0\times$\\
    matrix-free interleaved & --- & --- & 0.073 & 45.78 & 0.96 & $48\times $ & $4.6\times$\\
  \hline
  \end{tabular}
  \caption{Total solution time for different implementations of the PCG algorithm. Costs for matrix setup and host-device data transfer are listed separately and included in the total times. 100 iterations of the PCG main loop were carried out in all cases and all times are given in seconds.}
  \label{tab:TotalTimes}
 \end{center}
\end{table*}
\subsubsection{Time per iteration}\label{sec:TimePerIteration}
\begin{figure}
 \begin{center}
   \includegraphics[width=0.9\linewidth]{\figdir/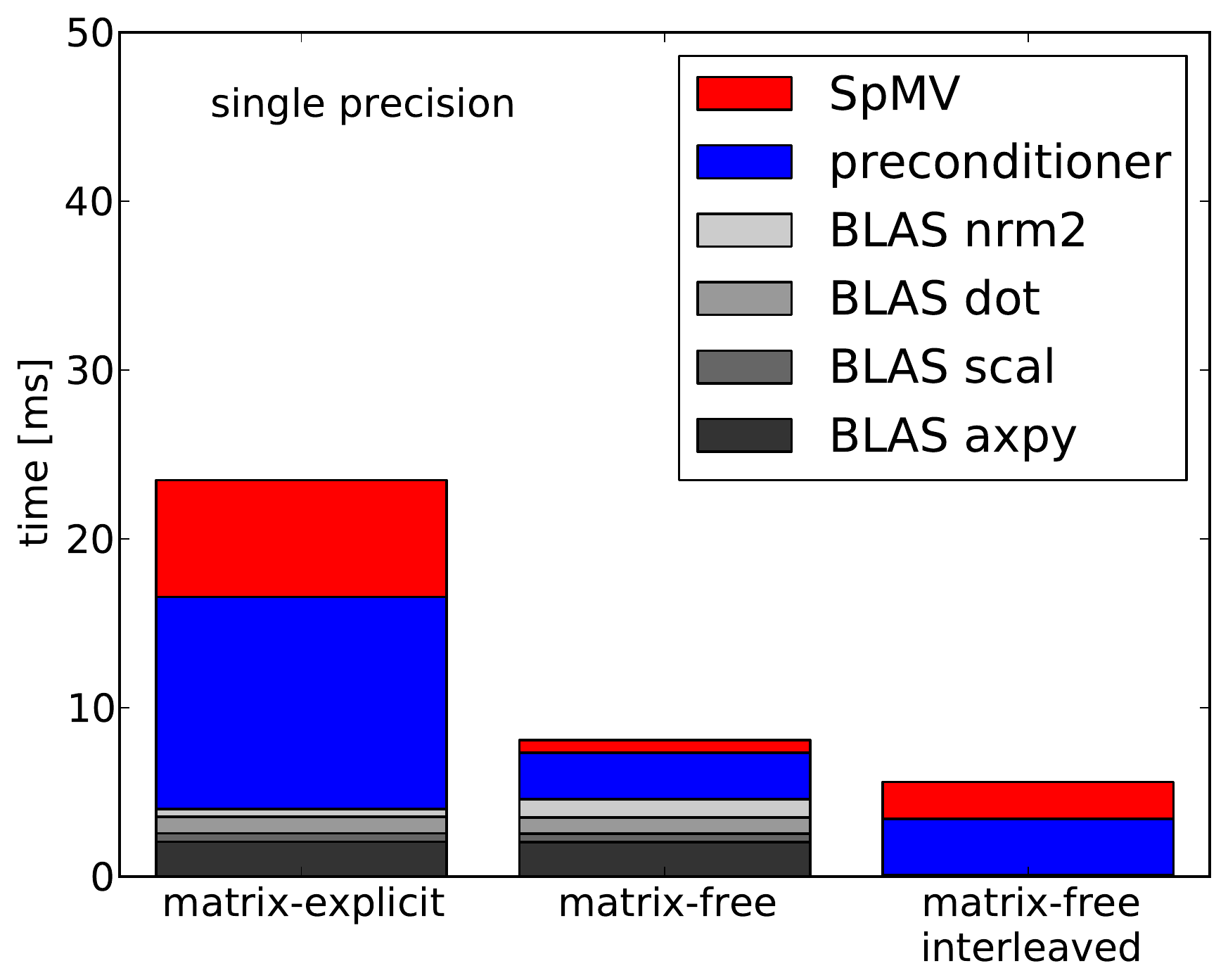}
 \end{center}
  \caption{Time per iteration for different parts of the main CG loop on the device using single precision. The BLAS-operations were implemented with the \CUBLAS\ library.}
  \label{fig:ResultsBreakdownCGLoopSingle}
\end{figure}

\begin{figure}
 \begin{center}
   \includegraphics[width=0.9\linewidth]{\figdir/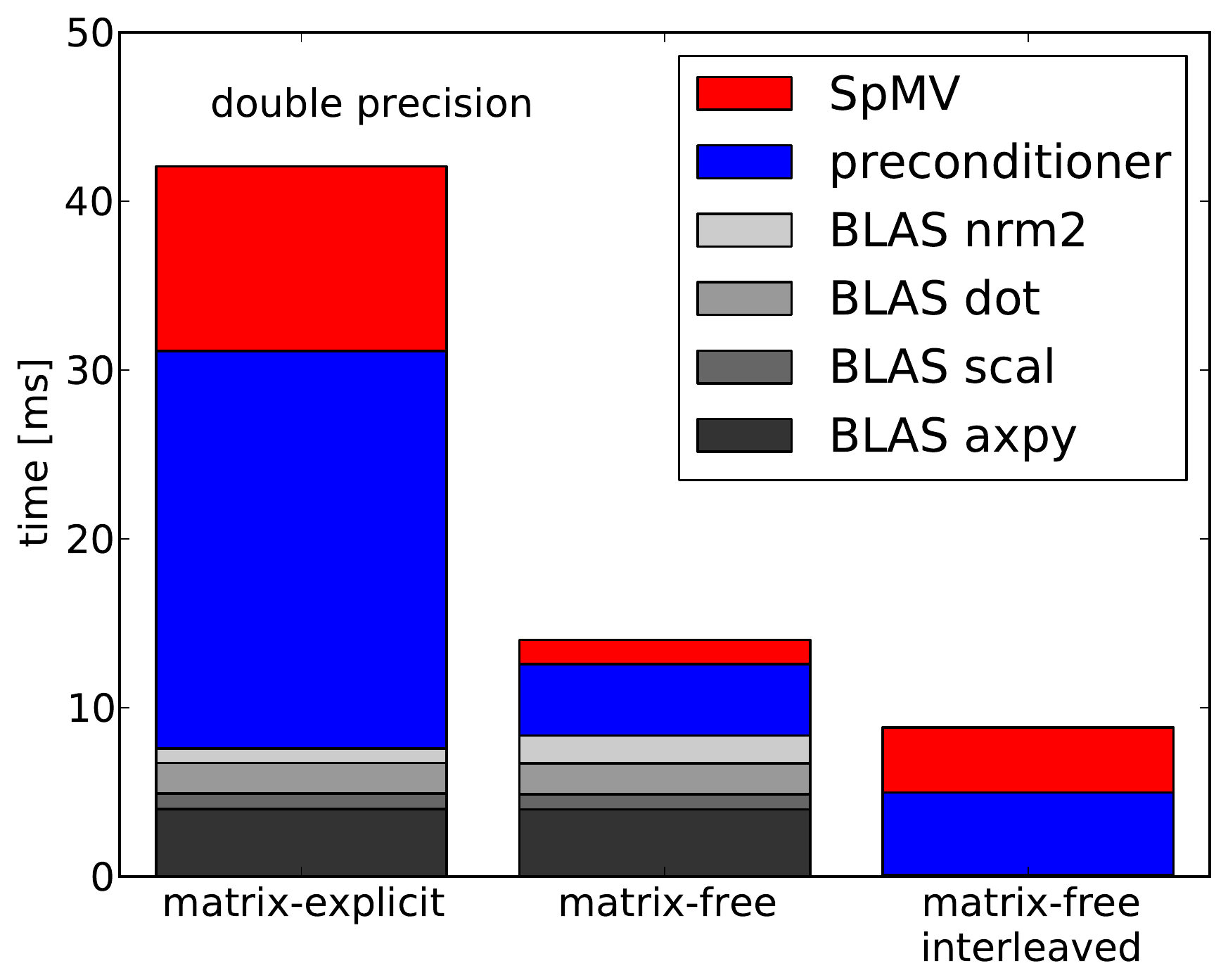}
 \end{center}
  \caption{Time per iteration for different parts of the main CG loop on the device using double precision. The BLAS-operations were implemented with the \CUBLAS\ library.}
  \label{fig:ResultsBreakdownCGLoopDouble}
\end{figure}

The time per iteration is given for the three different implementations of the PCG algorithm in Tab. \ref{tab:TimePerIteration} where we also calculate the speedup of the matrix-free implementation relative to the matrix-explicit version.
The speedup relative to the C implementation is $66\times$ for the matrix-free interleaved implementation in single precision and $52\times$ in double precision. In both cases it is more than four times faster than the matrix-explicit implementation on the GPU.
Comparing the total time per iteration for the matrix-free and the interleaved implementation, we find that the ratio between the two is $0.69$ in single precision and $0.63$ in double precision, which should be compared to the corresponding ratio of memory references in Tab. \ref{tab:PCGOperations}, which is $33/40=0.83$ without caching or $20/27=0.74$ assuming that the vectors $\vec{a}'$, $\vec{b}'$, $\vec{c}'$ and $\vec{d}$ and all data in one vertical column is cached (the ratio is $25/32=0.78$ if only the $\vec{a}'$, $\vec{b}'$, $\vec{c}'$ and $\vec{d}$ are cached). 
\begin{table}
 \begin{center}
  \begin{tabular}{|l|cc|cc|}
   \hline
   & \multicolumn{2}{|c|}{time per iteration} & \multicolumn{2}{|c|}{speedup}\\
             &   &      & C vs. & mat.--free \\
            implementation & C & CUDA &  CUDA & vs. CSR \\
\hline\hline
   \textit{single precision} & & & & \\
   matrix-explicit & 410.8 & 23.5 & $17\times$ & --- \\
   matrix-free & 376.3 & 8.1 & $46\times$ & $2.9\times$\\
   ---\textquotedbl--- interlvd. & 370.6 & 5.6 & $66\times$ & $4.2\times$\\
  \hline
   \textit{double precision} & & & & \\
   matrix-explicit & 494.6 & 41.4 & $12\times$ & --- \\
   matrix-free & 491.5 & 13.9 & $35\times$ & $3.0\times$\\
   ---\textquotedbl--- interlvd. & 453.4 & 8.8 & $52\times$ & $4.7\times$\\
  \hline
  \end{tabular}
  \caption{Time per iteration for different implementations of the conjugate gradient algorithm. All times are measured in milliseconds.}
  \label{tab:TimePerIteration}
 \end{center}
\end{table}
To identify the remaining bottlenecks, these times are further broken down for the GPU code in Figs. \ref{fig:ResultsBreakdownCGLoopSingle} and \ref{fig:ResultsBreakdownCGLoopDouble} for single- and double precision arithmetic. Note that for the matrix-free code the BLAS operations, and in particular the \verb!axpy! updates, make up a significant proportion of the time in the matrix-free code. The plot clearly shows the benefit of the interleaved implementation, which increases the performance by $28\%$ in single precision and $35\%$ in double precision.
\subsection{Absolute performance}
While comparing the runtime of the CUDA-C code to the corresponding sequential CPU implementation can give an idea of the achievable performance gains, it is somewhat arbitrary in that it depends on the exact CPU which is used for this comparison. For this reason we also quantified the absolute performance of the matrix-free CUDA-C code.

Several performance indicators for the kernels of the matrix-free code are shown in Tab. \ref{tab:ResultsPerformancePrecMatrixFree}. The global load efficiency measures the amount of coalesced global memory access and the L1 hit rate quantifies the cache efficiency. For all kernels load efficiency is very high due to coalescence of global memory access as described in section \ref{sec:MemoryLayout}.

\begin{figure}
 \begin{center}
   \includegraphics[width=0.9\linewidth]{\figdir/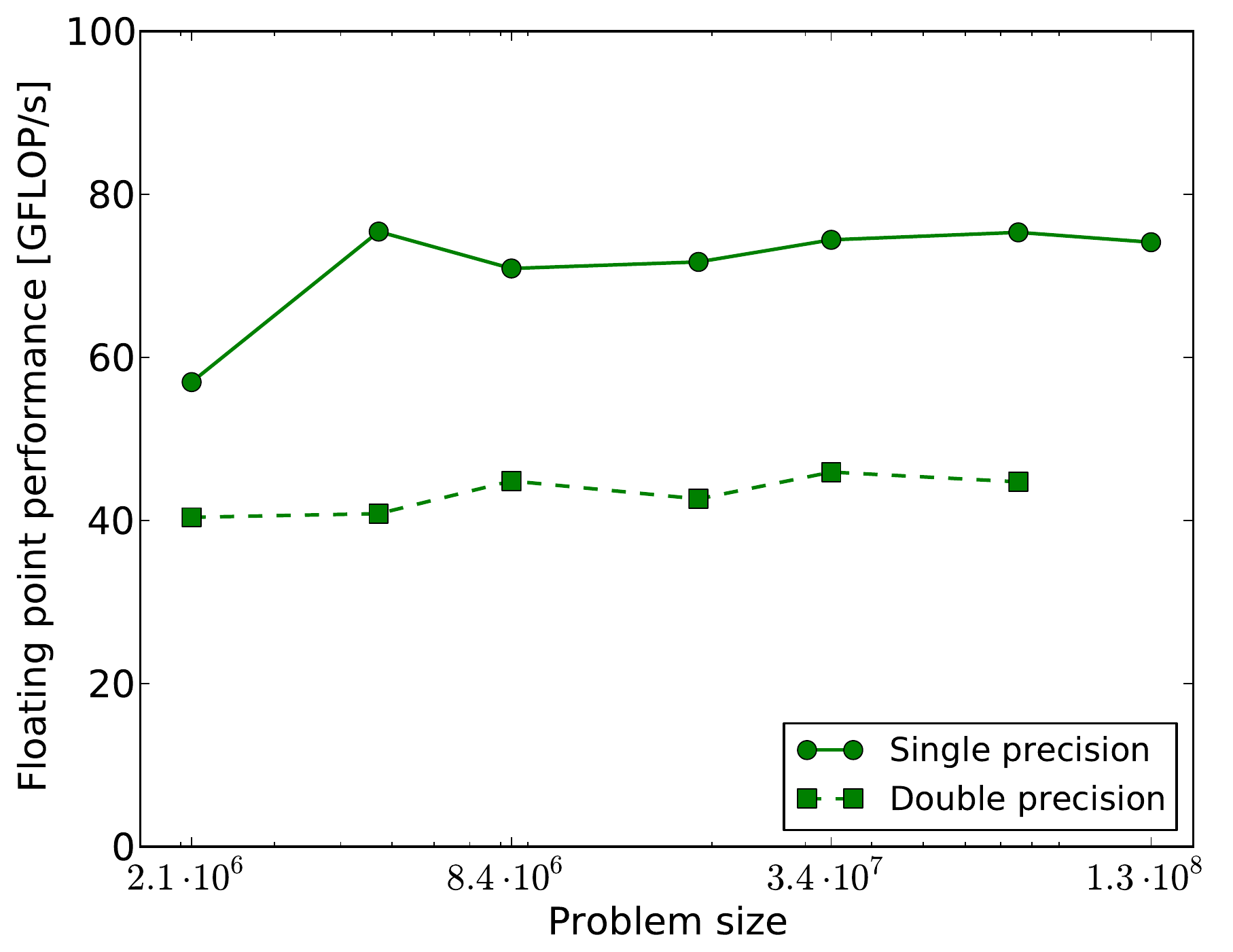}
 \end{center}
  \caption{Floating point performance for different problem sizes for the matrix-free interleaved PCG code on the GPU.}
  \label{fig:ResultsFloatingPointPerformance}
\end{figure}

\begin{figure}
 \begin{center}
   \includegraphics[width=0.9\linewidth]{\figdir/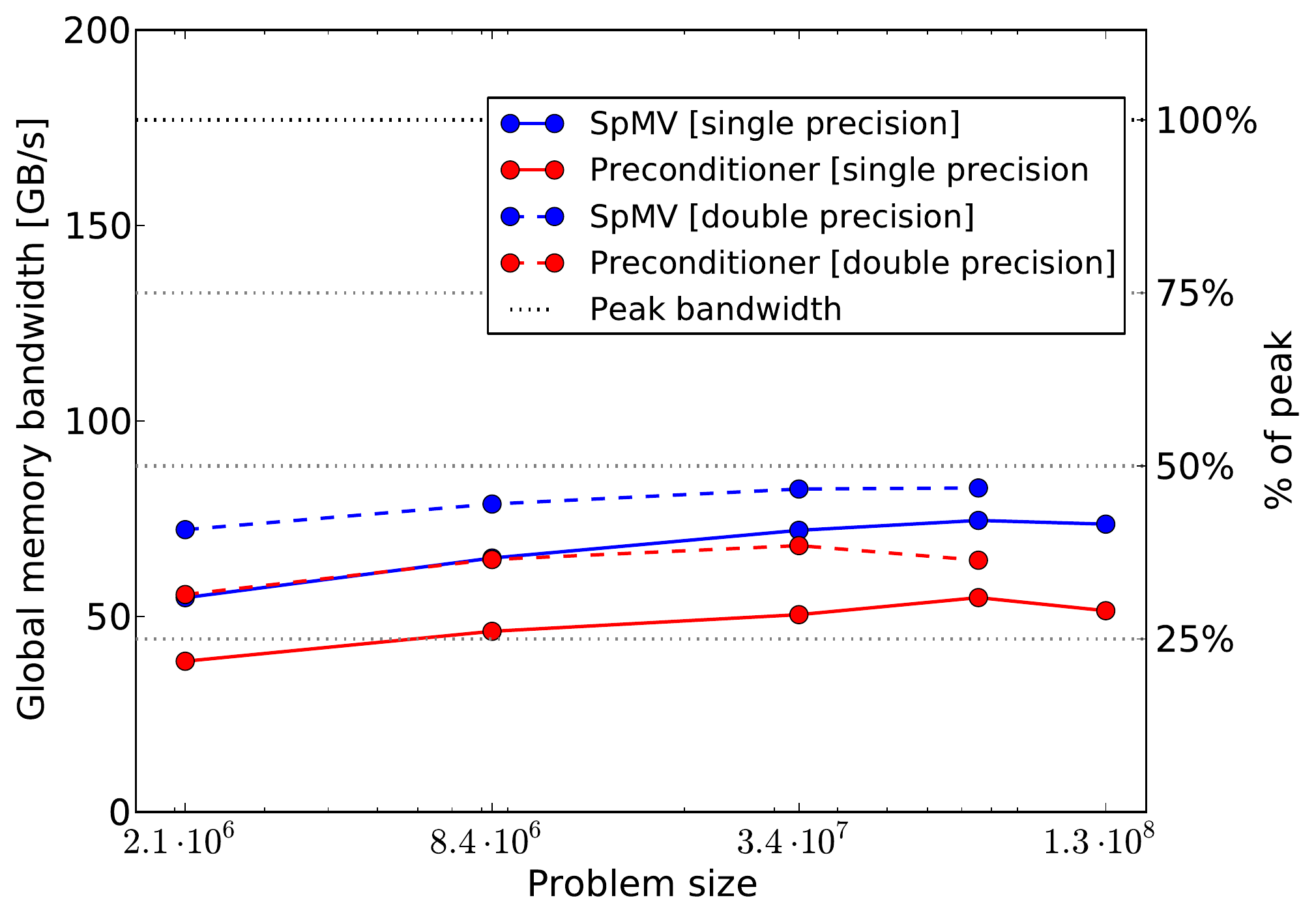}
 \end{center}
  \caption{Global memory bandwidth (DRAM read bandwidth as reported by the nVidia Visual Profiler) for different problem sizes for the interleaved apply- and preconditioner kernels.}
  \label{fig:ResultsBandwidth}
\end{figure}

The floating point performance for one iteration of the matrix-free interleaved PCG algorithm is plotted for a range of problem sizes between $2.1\cdot 10^6$ and $1.3\cdot10^8$ in Fig. \ref{fig:ResultsFloatingPointPerformance}. The nVidia M2090 Fermi GPU has a global memory of 6GB, and as the PCG algorithm requires the storage of 5 field vectors, which limits the problem size to less than $3\cdot10^8$ for single precision and $1.5\cdot10^8$ for double precision. In all cases the size of a vertical column was kept fixed at $n_z=128$. The performance is virtually independent of the problem size and about 70-80 GFLOPs for single- and 40-50 GFLOPs for double precision. As the algorithm is memory bound, a more relevant measure is the global memory band width which is shown for the interleaved kernels in Fig. \ref{fig:ResultsBandwidth}. The bandwidth increases slightly with the problem size and $25\%-50\%$ of the peak value could be achieved.
\begin{table}
 \begin{center}
  \begin{tabular}{|lrrrr|}
  \hline
    & GFLOPs & \multicolumn{1}{l}{memory}& \multicolumn{1}{l}{global} & \multicolumn{1}{l|}{L1} \\
    &        & \multicolumn{1}{l}{bandwidth} & \multicolumn{1}{l}{load ef-} & \multicolumn{1}{l|}{hit}\\
    &        & \multicolumn{1}{l}{[GB/s]}& \multicolumn{1}{l}{ficency} & \multicolumn{1}{l|}{rate}\\
\hline\hline
    \textit{single precision} & & & & \\
    SpMV         & 223.7 & 50.73 & 86.2\% & 76.7\% \\
    preconditioner & 45.4 & 38.41 & 99.9\% & 40.4\% \\
    interlvd. SpMV & 108.7 & 64.67 & 88.8\% & 56.5\% \\
    interlvd. prec. & 47.7 & 46.16 & 99.6\% & 33.6\% \\
\hline
    \textit{double precision} &&&&\\
    SpMV         & 119.0 & 79.07 & 85\% & 74.0\% \\
    preconditioner & 29.1 & 49.35 & 99.7\% & 40.1\% \\
    interlvd. SpMV & 60.9 & 78.97 & 88.0\% & 57.5\% \\
    interlvd. prec. & 32.8 & 64.58 & 99.8\% & 33.1\% \\
  \hline
  \end{tabular}
  \caption{Performance measurements of the different matrix-free kernels as reported by the nVidia visual profiler. The bandwidth is the DRAM load bandwidth.}
  \label{tab:ResultsPerformancePrecMatrixFree}
 \end{center}
\end{table}
\section{Conclusions}\label{sec:Conclusion}
In this article we described the matrix-free implementation of a Preconditioned Conjugate Gradient solver for strongly anisotropic elliptic PDEs on a GPU. 
Equations of this type arise in many applications in atmospheric and geophysical modelling if the problem is discretised in ``flat'' geometries. In particular the semi-implicit semi\hyp Lagrangian time discretisation of the non-hydrostatic Euler equations in the dynamical core of many weather- and climate prediction models leads to a three dimensional PDE for the pressure correction and due to the small thickness of the atmosphere the elliptic operator has a very strong anisotropy in the vertical direction. This anisotropy can be exploited to construct a simple and efficient preconditioner based on vertical line relaxation. The dependencies in each vertical column require a horizontal domain decomposition, which has implementations for the data- and thread layout on the GPU.

As the performance of the algorithm is limited by the speed with which data can transferred from global memory to the compute units, it is important to reduce the number of memory references. We achieved this by using a matrix-free implementation which recalculates the local matrix stencil whenever it is needed instead of reading it from memory. In addition, we reduced the amount of data transfer by fusing several kernels in the PCG loop, which gave an additional improvement of $28\%$ in single precision and $35\%$ in double precision. In total we demonstrated that on an nVidia Fermi M2090 GPU the best matrix-free code achieved a speedup of $60\times$ in single precision and $48\times$ for double precision, compared to the corresponding sequential implementation on an Intel Xeon E5-2620 Sandybridge CPU. In terms of absolute times, the residual for a problem with $8.3\cdot 10^6$ degrees of freedom could be reduced by more than five orders of magnitude in 0.62 seconds in single precision. The double precision implementation takes about 1.5 times longer, which demonstrates the good double precision performance of modern GPUs. Our matrix-free version is more than four times faster than a matrix-explicit GPU implementation based on the \CUSPARSE\ and \CUBLAS\ libraries using the CSR matrix format which clearly demonstrates the benefit of our approach. We measured the absolute performance of our code for a range of problem sizes and achieved a global memory bandwidth of $25\%-50\%$ of the peak rate for problem sizes between $2.1\cdot10^6$ and $1.3\cdot10^8$ degrees of freedom. 

Global memory access was coalesced for the threads within a half warp by numbering the degrees of freedom such that the horizontal index runs fastest. This is different from CPU implementations where vertical columns are stored consecutively in memory for cache efficiency. While the specific implementation discussed in this article is based on a regular horizontal grid, our method can be applied to any three dimensional grid which can be written as the tensor product of possibly unstructured horizontal grid and a non-uniform one dimensional grid in the vertical direction.

The achieved global memory bandwidth is a sizeable fraction of the peak value, and it can be theoretically improved by an additional factor $2\times$ to $4\times$ by making better use of the GPU cache or shared memory. This would require the parallelisation of the tridiagonal solver in the vertical direction, and we are currently investigating the substructuring approach described in \cite{Toselli2005}.

With the planned increase in weather- and climate model resolution, problems with more than $10^{10}$ degrees of freedom need to be solved. Clearly this is not possible on a single GPU due to limited global memory size. To solve problems of this size hundreds of GPUs are necessary as each has limited memory. We are currently extending the algorithm to multi-GPU clusters will introduce additional bottlenecks: unless data can be copied directly between GPU memory, it has to be transferred to the host at each iteration before it can be sent through the standard MPI network. However, in this case only halo data needs to be exchanged between neighbouring devices and given that the local problem size is not too small, this is significantly less than the total data processed by one GPU.

The PCG solver described in this work requires around hundred iterations to reduce the residual by five orders of magnitude. In contrast, multigrid solvers can achieve the same reduction in a much smaller number of iterations, as has been demonstrated in \cite{Mueller2013} for the elliptic PDE considered in this article. The preconditioner used in this work can be used as a multigrid smoother and the only missing components are intergrid operators for restriction and prolongation, which is the object of our current research.
\ifpreprint 
\section*{Acknowledgements}
\else 
\begin{acknowledgements}
\fi 
We would like to thank all members of the GungHo! project and in particular Chris Maynard and David Ham for useful and inspiring discussions. The numerical experiments in this work were carried out on a node of the aquila supercomputer at the University of Bath and we are grateful to Steven Chapman for his continuous and tireless technical support which was essential for the success of this project. The contribution of EM and RS was funded as part of the NERC project on ``Next Generation Weather and Climate Prediction'' (NGWCP), grant number NE/J005576/1. 
\ifpreprint 
\else 
\end{acknowledgements}
\fi 
\appendix
\section{Interleaved PCG kernels}\label{sec:Algorithms}
The kernels for the interleaved PCG algorithm described in section \ref{sec:PCGilvd} are shown in Algorithms \ref{alg:InterleavedApply} and \ref{alg:InterleavedPrec}. In the GPU implementation each thread calculates dot products and norms in one column. To obtain global sums, these need to be reduced with an additional BLAS operation. However, as this operation only operates on two-dimensional (horizontal) vectors, its cost is very small ($<1\%$ of the time per iteration).
\begin{algorithm}
  \caption{Interleaved matrix-multiplication kernel.
      Simultaneously calculate $\vec{u}\mapsfrom \vec{u}+\alpha \vec{p}$, $\vec{p}\mapsfrom \vec{z}+\beta \vec{p}$, $\vec{q} \mapsfrom \matrix{A}\vec{z}+\beta \vec{q}$, $\sigma\mapsfrom\langle \vec{p},\vec{q}\rangle$
      in a single iteration over the grid.}
  \label{alg:InterleavedApply}
  \begin{algorithmic}[1]
    \FOR{$i=0,\dots,m$}
    \FOR{$j=0,\dots,m$}
    \STATE{Calculate $\alpha_{i',j'}$ and $|T_{ij}|$}
    \STATE{$\sigma\mapsfrom 0$}
    \FOR{$k=0,\dots,n_z-1$}
       \STATE{$p^*\mapsfrom p_{ijk}$,
              $q^*\mapsfrom q_{ijk}$, 
              $z^*\mapsfrom z_{ijk}$}
       \STATE{$u_{ijk}\mapsfrom u_{ijk}+\alpha\cdot p^*$}
       \STATE{$p^*\mapsfrom \beta\cdot p^*+z^*$,
              $q^*\mapsfrom \beta\cdot q^*$}
       \STATE{$p_{ijk}\mapsfrom p^*$}
       \STATE{$\delta q\mapsfrom \left((a'_k-b'_k-c'_k)\cdot|T_{ij}|-\alpha_{ij}\right)\cdot z^*$\\
              $\quad+\;b'_k\cdot |T_{ij}|\cdot z_{i,j,k+1}+c'_k\cdot |T_{ij}|\cdot z_{i,j,k-1}$\\
              $\quad+\;\alpha_{i+1,j}\cdot z_{i+1,j,k}+\alpha_{i-1,j}\cdot z_{i-1,j,k}$\\$\quad+\;\alpha_{i,j+1}\cdot z_{i,j+1,k}+\alpha_{i,j-1}\cdot z_{i,j-1,k}$}
       \STATE{$q^*\mapsfrom q^*+d_k\cdot \delta q$,
              $\sigma\mapsfrom \sigma+p^*\cdot q^*$}
       \STATE{$q_{ijk}\mapsfrom q^*$}
    \ENDFOR
    \ENDFOR
    \ENDFOR
  \end{algorithmic}
\end{algorithm}
\begin{algorithm}
  \caption{Interleaved preconditioner kernel. Simultaneously calculate
      $\vec{r}\mapsfrom \vec{r}-\alpha \vec{q}$,
      $R\equiv||\vec{r}||\mapsfrom \sqrt{\langle \vec{r},\vec{r}\rangle}$,
	    $\kappa\mapsfrom  \langle \vec{r},\vec{z}\rangle$ and solve 	    $\matrix{M}\vec{z}=\vec{r}$ in a single iteration over the grid.}
  \label{alg:InterleavedPrec}
  \begin{algorithmic}[1]
    \FOR{$i=0,\dots,m$}
    \FOR{$j=0,\dots,m$}
    \STATE{Calculate $\alpha_{i',j'}$ and $|T_{ij}|$}
    \STATE{$R\mapsfrom 0, \kappa\mapsfrom 0$,
           $D\mapsfrom (a'_k-b'_k-c'_k)-\alpha'_{ij}$,
           $\phi_0\mapsfrom b'_k/D$}
    \STATE{$r^*\mapsfrom r_{ij0}-\alpha\cdot q_{ij0}$,
           $R\mapsfrom R+r^*\cdot r^*$}
    \STATE{$z_{ij0}\mapsfrom r^*/(D\cdot |T_{ij}|\cdot d_k)$,
           $r_{ij0}\mapsfrom r^*$}
    \FOR{$k=0,\dots,n_z-1$}
      \STATE{$D\mapsfrom \left((a'_k-b'_k-c'_k)-\alpha'_{ij}\right)-\phi_{k-1}\cdot c'_k$,
             $\phi_k\mapsfrom b'_k/D$}
      \STATE{$r^*\mapsfrom r_{ijk}-\alpha\cdot q_{ijk}$,
             $R\mapsfrom R+r^*\cdot r^*$}
      \STATE{$z_{ijk}\mapsfrom \left(r^*/(|T_{ij}|\cdot d_k)-c'_k\cdot z_{i,j,k-1}\right)/D$,
             $r_{ijk}\mapsfrom r^*$}
    \ENDFOR
    \STATE{$\kappa\mapsfrom \kappa + z_{i,j,n_z-1}\cdot r_{i,j,n_z-1}$}
    \FOR{$k=n_z-2,\dots,0$}
      \STATE{$z^*\mapsfrom z_{ijk}-\phi_k\cdot z_{i,j,k+1}$,
             $\kappa\mapsfrom \kappa+z^*\cdot r_{ijk}$}
      \STATE{$z_{ijk}\mapsfrom z^*$}
    \ENDFOR
    \ENDFOR
    \ENDFOR
  \STATE{$R\mapsfrom \sqrt{R}$}
\end{algorithmic}
\end{algorithm}
\ifpreprint 
\bibliographystyle{unsrt}
\else 
\bibliographystyle{spbasic}
\fi 

\end{document}